\documentclass[11pt]{article}
\usepackage{amsfonts}
\usepackage{color}
\newtheorem{theo}{Theorem}[section]
\newtheorem{lem}[theo]{Lemma}

\newcommand{\mysection}[1]{\section{#1} \setcounter{equation}{0}}
\newcommand{\proof}{{\sc Proof.} \quad}
\newcommand{\proofc}{{\sc Proof} \ }
\newcommand{\be}{\begin{equation} \label}
\newcommand{\ee}{\end{equation}}
\newcommand{\bea}{\begin{eqnarray}\label}
\newcommand{\eea}{\end{eqnarray}}
\newcommand{\bas}{\begin{eqnarray*}}
\newcommand{\eas}{\end{eqnarray*}}
\newcommand{\bit}{\begin{itemize}}
\newcommand{\eit}{\end{itemize}}
\newcommand{\qed}{\hfill$\Box$ \vskip.2cm}
\newcommand{\nn}{\nonumber}
\newcommand{\R}{\mathbb{R}}
\newcommand{\N}{\mathbb{N}}
\newcommand{\pO}{\partial\Omega}

\newcommand{\eps}{\varepsilon}

\newcommand{\abs}{\\[5pt]}

\newcommand{\cb}{\color{blue}}

\newcommand{\io}{\int_\Omega}
\newcommand{\tm}{T_{max}}

\newcommand{\uw}{\underline{w}}
\newcommand{\ow}{\overline{w}}

\newcommand{\wout}{w_{out}}
\newcommand{\win}{w_{in}}
\newcommand{\parab}{{\cal P}}

\newcommand{\al}{a_\lambda}
\newcommand{\bl}{b_\lambda}
\newcommand{\psil}{\psi_\lambda}

\def\va{\raise 2pt\hbox{,}}
\begin{document}
\title{A degenerate chemotaxis system with flux limitation:\\
 Finite-time blow-up}
\author{
Nicola Bellomo\footnote{nicola.bellomo@polito.it}\\
{\small Department of Mathematics, Faculty  Sciences,}\\
 {\small King Abdulaziz University,  Jeddah, Saudi Arabia}\\
      {\small Politecnico di Torino,  10129 Torino, Italy}\\
\and
Michael Winkler\footnote{michael.winkler@math.uni-paderborn.de}\\
{\small Institut f\"ur Mathematik, Universit\"at Paderborn,}\\
{\small 33098 Paderborn, Germany} }
\date{}
\maketitle
\begin{abstract}
\noindent
  This paper is concerned with radially symmetric solutions of the parabolic-elliptic version
  of the Keller-Segel system with flux limitation, as given by
  \bas
	\left\{ \begin{array}{l}
\displaystyle u_t=\nabla \cdot \Big(\frac{u\nabla u}{\sqrt{u^2+|\nabla u|^2}}\Big) - \chi \, \nabla \cdot \Big(\frac{u\nabla v}{\sqrt{1+|\nabla v|^2}}\Big), \\[1mm]
	0=\Delta v - \mu + u,
	\end{array} \right.
	\qquad \qquad (\star)
  \eas
  under the initial condition $u|_{t=0}=u_0>0$ and no-flux boundary conditions
  in a ball $\Omega\subset\R^n$, where $\chi>0$ and $\mu:=\frac{1}{|\Omega|} \io u_0$.
  A previous result\cite{bellomo-winkler1}
  has asserted global existence of bounded classical solutions for arbitrary
  positive radial initial data $u_0\in C^3(\bar\Omega)$ when
  either $n\ge 2$ and $\chi<1$, or $n=1$ and $\io u_0<\frac{1}{\sqrt{(\chi^2-1)_+}}$.\abs
  This present paper shows that these conditions are essentially optimal: Indeed,
  it is shown that if the taxis coefficient is large enough in the sense that $\chi>1$,
  then for any choice of
  \bas
	\left\{ \begin{array}{ll}
	m>\frac{1}{\sqrt{\chi^2-1}} \quad & \mbox{if } n=1, \\[2mm]
	m>0 \mbox{ is arbitrary } \quad & \mbox{if } n\ge 2,
	\end{array} \right.
  \eas
  there exist positive initial data $u_0\in C^3(\bar\Omega)$ satisfying $\io u_0=m$ which are such that
  for some $T>0$,
  ($\star$) possesses a uniquely determined classical solution $(u,v)$ in $\Omega\times (0,T)$ blowing up
  at time $T$ in the sense that $\limsup_{t\nearrow T} \|u(\cdot,t)\|_{L^\infty(\Omega)}=\infty$.\abs
  This result is derived by means of a comparison argument applied to the doubly degenerate scalar parabolic
  equation satisfied by the mass accumulation function associated with ($\star$).\abs
\end{abstract}

\newpage
\mysection{Introduction}
{\bf Flux-limited Keller-Segel systems.} \quad
This paper presents a continuation of the analytical study \cite{bellomo-winkler1}
of a flux limited chemotaxis model recently derived as a development of
the classical pattern formation model proposed by Keller and Segel (\cite{keller_segel}) to model collective behavior
of populations mediated by a chemoattractant.
In a general form, this model describes the spatio-temporal evolution of the cell density $u=u(x,t)$ and
the chemoattractant concentration $v=v(x,t)$ by means of the parabolic system
%
%he general structure of this new model is as follows:
%, here in after the (Flim-KS) model, is as follows:
\begin{equation}\label{Flim-KS}
\left\{
\begin{array}{l}
\displaystyle u_t =  \nabla \cdot \left(D_u (u,v)\frac{u \nabla u}{\sqrt{u^2 + |\nabla u|^2}} -  S(u,v)
\frac{u \nabla v}{\sqrt{1+|\nabla {v}|^2}}\right) +H_1(u,v)\va\\
 \\
\displaystyle v_t =  D_v \Delta v + H_2(u,v),
\end{array} \right.
\end{equation}
where $D_u$ and $D_v$ denote the respective diffusivities, $S$ represents the chemotactic sensitivity and
$H_1$ and $H_2$ account for mechnisms of proliferation, degradation, and possibly also
interaction.
In comparison to the original Keller-Segel system,
besides including cell diffusivities inhibited at small densities and hence supporting finite propagation speeds,
the main innovative aspect in (\ref{Flim-KS}) apparently
consists in the choice of limited diffusive and cross-diffusive fluxes 	%and of the optimal transport \cite{Brenier-09}
in the first equation.
From a modeling perspective, this is founded in the assumption that particles do not diffuse
arbitrarily in the space but, on the contrary, through some privileged ways such as the border of cells,
and the particular mathematical structure in (\ref{Flim-KS}) has been derived in \cite{bellomo-bellouquid-nieto-soler2010}
by asymptotic limits and time-space scaling on the description delivered by kinetic-type models, where
cell-cell interactions are modeled by theoretical tools of evolutive game theory.
(see also \cite{BBTW}; cf.~the survey \cite{hillen-painter-09}, and also \cite{Horstmann-03},
for a review of modeling issues based on the classical approach of continuum mechanics).\abs
{\bf Blow-up in semilinear and quasilinear chemotaxis systems.} \quad
The goal of the present work is to clarify to what extent the introduction of such flux limitations may suppress
phenomena of blow-up, as known to constitute one of the most striking characteristic features of the classical
Keller-Segel system
\be{KS}
	\left\{ \begin{array}{l}
	u_t=\Delta u - \nabla \cdot (u\nabla v), \\[1mm]
	v_t=\Delta v - v + u,
	\end{array} \right.
\ee
and also of several among its derivates.
Indeed, the Neumann initial-boundary value problem for (\ref{KS})
is known to possess solutions blowing up in finite time with respect to the spatial $L^\infty$
norm of $u$ when either the spatial dimension $n$ satisfies
$n\ge 3$ (\cite{win_JMPA}), or when $n=2$ and the initially present -- and thereafter conserved --
total mass $\int u(\cdot,0)$ of cells is suitably large
(\cite{herrero_velazquez}, \cite{mizoguchi_win}).
On the other hand, in the case $n=2$ appropriately small values of $\int u(\cdot,0)$
warrant global existence of bounded solutions \cite{nagai_senba_yoshida},
whereas if $n\ge 3$ then global bounded solutions exist under alternative smallness conditions involving the norms of
$(u(\cdot,0),v(\cdot,0))$ in $L^\frac{n}{2} \times W^{1,n}$ (\cite{cao_optspace}, \cite{win_JDE2010}).
In the associated spatially one-dimensional problem, global bounded solutions exist for all reasonably regular
initial data, thus reflecting absence of any blow-up phenomenon in this case (\cite{osaki_yagi}).\abs
The knowledge on corresponding features of quasilinear relatives of (\ref{KS}) seems most developed for models
involving density-dependent variants in the diffusivity and the chemotactic sensitivity.
For instance, if $D_u$ and $S$ are smooth positive functions on $[0,\infty)$, then the Neumann problem for
\be{KSn}
	\left\{ \begin{array}{l}
	u_t=\nabla \cdot (D_u(u)\nabla u) - \nabla \cdot (S(u)\nabla v), \\[1mm]
	v_t=\Delta v - v + u,
	\end{array} \right.
\ee
posesses some unbounded solutions whenever $\frac{S(u)}{D_u(u)} \ge Cu^{\frac{2}{n}+\eps}$ for all $u\ge 1$ and some
$C>0$ and $\eps>0$ (\cite{win_volume});
beyond this,
refined studies have given additional conditions on $D_u$ and $S$ under which this
singularity formation must occur within finite time, and have moreover
identified some particular cases of essentially algebraic behavior of both $D_u$ and $S$ in which these explosions
must occur in infinite time only (\cite{cieslak_stinner_JDE2012}, \cite{cieslak_stinner_ACAP},
\cite{cieslak_stinner_JDE2015}, see also
\cite{cieslak_laurencot_ANNIHP} for a related example on finite-time blow-up).
The optimality of the above growth condition is indicated by a result in \cite{taowin_JDE_subcrit} and
\cite{ishida_seki_yokota_JDE2014} asserting global
existence of bounded solutions in the case when $\frac{S(u)}{D_u(u)} \le Cu^{\frac{2}{n}-\eps}$ for $u\ge 1$ with some
$C>0$ and $\eps>0$, provided that $D_u$ decays at most algebraically as $u\to\infty$ (cf.~e.g.~\cite{kowalczyk_szymanska},
\cite{horstmann_win}, \cite{wrzosek} and \cite{senba_suzuki_AAA}
for some among the numerous precedents in this direction).\abs
As compared to this, the literature on variants of (\ref{KS})
involving modifications of the dependence of fluxes on gradients seems quite thin.
Moreover, the few results available in this direction mainly seem to concentrate on modifications in the
cross-diffusive term, essentially guided by the underlying idea
to rule out blow-up by suitable regularizations of the taxis term in (\ref{KS}),
as apparently justified in appropriate biological contexts (see the discussion in \cite{hillen-painter-09}
as well as the analytical findings reported there).
In particular, we are not aware of any result detecting an explosion in any such context;
this may reflect the evident challenges connected to rigorously proving the occurrence of blow-up in
such complex chemotaxis systems.\abs
{\bf Main results: Detecting blow-up under optimal conditions.} \quad
The present work will reveal that actually also the introduction of flux limitations need not necessarily
suppress phenomena of chemotactic collapse in the sense of blow-up.
In order to make this manifest in a particular setting, let us concentrate on the case when in (\ref{Flim-KS})
we have $D_u\equiv 1$ and $S\equiv const.$ as well as $H_1\equiv 0$, and in order to simplify our analysis let us
moreover pass to a parabolic-elliptic simplification thereof, thus focusing on a frequently considered limit case of fast
signal diffusion (\cite{jaeger_luckhaus}).
Here we note that e.g.~in the previously discussed situations of (\ref{KS}) and (\ref{KSn}), up to
few exceptions (\cite{biler_corrias_dolbeault}) such parabolic-elliptic
variants are known to essentially share the same properties as the respective fully parabolic model
with regard to the occurrence of blow-up (\cite{nagai2001}, \cite{biler1999}, \cite{calvez_carrillo}, \cite{djie_win}).\abs
We shall thus subsequently be concerned with the initial-boundary value problem
\be{0}
	\left\{ \begin{array}{l}
	u_t=\nabla \cdot \Big(\frac{u\nabla u}{\sqrt{u^2+|\nabla u|^2}}\Big)
		- \chi \nabla \cdot \Big(\frac{u\nabla v}{\sqrt{1+|\nabla v|^2}}\Big),
	\qquad x\in \Omega, \ t>0, \\[1mm]
	0=\Delta v - \mu + u,
	\qquad x\in \Omega, \ t>0, \\[1mm]
	\Big(\frac{u\nabla u}{\sqrt{u^2+|\nabla u|^2}}
		- \chi \frac{u\nabla v}{\sqrt{1+|\nabla v|^2}} \Big) \cdot \nu=0,
	\qquad x\in \partial\Omega, \ t>0, \\[1mm]
	u(x,0)=u_0(x),
	\qquad x\in\Omega,
 	\end{array} \right.
\ee
in a ball $\Omega=B_R(0)\subset \R^n$, $n\ge 1$, where $\chi>0$ and the initial data are such that
\be{init}
	u_0\in C^3(\bar\Omega)
	\qquad \mbox{is radially symmetric and positive in $\bar\Omega$ with $\frac{\partial u_0}{\partial\nu}=0$ on } \pO,
\ee
and where
\be{mu}
	\mu:=\frac{1}{|\Omega|} \io u_0(x)dx
\ee
denotes the spatial average of the latter.\abs
In fact, it has been shown in \cite{bellomo-winkler1} that this problem is well-posed, locally in time,
in the following sense.\abs
%
%
%\begin{theo}\label{theo43}
{\bf Theorem A} \
{\it
  Let $n\ge 1, \chi>0$ and $\Omega:=B_R(0)\subset \R^n$ with some $R>0$, and
  suppose that $u_0$ complies with (\ref{init}). Then there exist $\tm \in (0,\infty]$ and a uniquely determined
  pair $(u,v)$ of positive radially symmetric functions
  $u\in C^{2,1}(\bar\Omega\times [0,\tm))$ and $v\in C^{2,0}(\bar\Omega\times [0,\tm))$
  which solve (\ref{0}) classically in $\Omega\times (0,\tm)$, and which are such that
  \be{43.1}
	\mbox{if $\tm<\infty$ \quad then \quad}
	\limsup_{t\nearrow\tm} \|u(\cdot,t)\|_{L^\infty(\Omega)}=\infty.
  \ee
}

%\end{theo}
%
%
%
%
Now in order to formulate our results and put them in perspective adequately,
let us moreover recall the following statement
on global existence and boundedness in certain subcritical cases
which has been achieved in \cite{bellomo-winkler1}.\abs
%
%
%\begin{theo}\label{theo555}
{\bf Theorem B} \
{\it
  Let $\Omega:=B_R(0)\subset \R^n$ with some $R>0$, and
  assume that $u_0$ satisfies (\ref{init}), and that either
  \be{555.01}
	n\ge 2
	\qquad \mbox{and} \qquad
	\chi<1,
  \ee
  or
  \be{555.02}
	n=1, \quad \chi>0
	\qquad \mbox{and} \qquad
	\io u_0<m_c,
  \ee
  where in the case $n=1$ we have set
  \be{mc}
	m_c:=\left\{ \begin{array}{ll}
	\frac{1}{\sqrt{\chi^2-1}} \qquad & \mbox{if } \chi>1, \\[3mm]
	+\infty & \mbox{if } \chi \le 1.
	\end{array} \right.
  \ee
%  Then for any choice of $u_0$ fulfilling (\ref{init}),
  Then the problem (\ref{0})
  possesses a unique global classical solution
  $(u,v)\in C^{2,1}(\bar\Omega\times [0,\infty)) \times C^{2,0}(\bar\Omega\times [0,\infty))$
  which is radially symmetric and such that for some $C>0$ we have
  \be{555.1}
	\|u(\cdot,t)\|_{L^\infty(\Omega)} \le C
	\quad \mbox{and} \quad
	\|v(\cdot,t)\|_{L^\infty(\Omega)} \le C
	\qquad \mbox{for all } t>0.
  \ee
}

%\end{theo}
%
%
%
%
It is the purpose of the present work to complement the above result on global existence by showing that in both cases
$n\ge 2$ and $n=1$, the conditions (\ref{555.01}) and (\ref{555.02}) are by no means artificial and of purely
technical nature, but that in fact they are essentially optimal in the sense
that if appropriate reverse inequalities hold, then finite-time blow-up may occur.
To be more precise, the main results of this paper can be formulated as follows.
\begin{theo}\label{theo14}
  Let $n\ge 1$ and $\Omega:=B_R(0)\subset \R^n$ with some $R>0$, and suppose that
  \be{chi}
	\chi>1,
  \ee
  and that
  \be{m}
	\left\{ \begin{array}{ll}
	m>m_c \quad & \mbox{if } n=1, \\[1mm]
	m>0 \mbox{ is arbitrary } \quad & \mbox{if } n\ge 2,
	\end{array} \right.
  \ee
  where $m_c$ is as in (\ref{mc}).
  Then there exists a nondecreasing function $M\in C^0([0,R])$ fulfilling $M(R)\le m$,
  which is such that whenever $u_0$ satisfies (\ref{init}) as well as
  \be{14.1}
	\int_{B_r(0)} u_0(x)dx \ge M(r)
	\qquad \mbox{for all } r\in [0,R],
  \ee
  the solution $(u,v)$ of (\ref{0}) blows up in finite time in the sense that in Theorem A	%\ref{theo43}
  we have $\tm<\infty$
  and
  \be{14.2}
	\limsup_{t\nearrow \tm} \|u(\cdot,t)\|_{L^\infty(\Omega)}=\infty.
  \ee
  In particular, for each prescribed number $m$ fulfilling (\ref{m}) one can find initial data $u_0$ which satisfy
  (\ref{init}) and $\io u_0=m$, and which are such that the corresponding solution of (\ref{0}) has the property (\ref{14.2}).
\end{theo}
In comparison to the classical Keller-Segel system (\ref{KS}), this in particular
means that when $n\ge 2$, the possible occurrence of blow-up does not go along with a critical mass phenomenon,
but that there rather exists a
{\em critical sensitivity parameter}, namely $\chi=1$, which distinguishes between existence and nonexistence
of blow-up solutions.
On the other hand, if $n=1$, then for any $\chi>1$, beyond this there exists a critical mass phenomenon,
in quite the same flavor as present in (\ref{KS}) when $n=2$.\abs
{\bf Plan of the paper.} \quad
Due to the apparent lack of an adequate global dissipative structure, a blow-up
analysis for (\ref{0}) cannot be built on the investigation of any energy functional,
as possible in both the original Keller-Segel system (\ref{KS}) and its quasilinear variant (\ref{KSn})
(\cite{win_JMPA}, \cite{mizoguchi_win}, \cite{cieslak_stinner_JDE2012}).
Apart from this, any reasoning in this direction needs to adequately cope with the circumstance that
as compared to (\ref{KS}), in (\ref{0}) the cross-diffusive flux is considerably inhibited wherever $|\nabla v|$ is large,
which seems to prevent access to blow-up arguments based on tracking the evolution of weighted $L^1$ norms of $u$
such as e.g.~the moment-like functionals considered in \cite{nagai2001}.\abs
That blow-up may occur despite this strong limitation of cross-diffusive flux will rather be shown by a comparison
argument. Indeed, it can readily be verified (Lemma \ref{lem_w})
that given a radial solution $u$ of (\ref{0}) in $B_R\times (0,T)$, the mass accumulation function $w=w(s,t)$,
as defined in a standard manner by introducing $w(s,t):=\int_0^{s^\frac{1}{n}} r^{n-1} u(r,t)dr$,
$(s,t)\in [0,R^n]\times [0,T)$,
satisfies a scalar parabolic equation which is doubly degenerate, both in space
as well as with respect to the variable $w_s$, but after all allows for an appropriate comparison principle
for certain generalized sub- and supersolutions (Lemma \ref{lem70}).\abs
Accordingly, at the core of our analysis will be the construction of suitable subsolutions to the respective problem;
in fact, we shall find such subsolutions $\uw$ which undergo a finite-time gradient blow-up at the origin in the sense that
for some $T>0$ we have $\sup_{s\in (0,R^n)}\frac{\uw(s,t)}{s} \to\infty$ as $t\nearrow T$, implying blow-up of $u$
before or at time $T$ whenever $w(\cdot,0)$ lies above $\uw(\cdot,0)$.
These subsolutions will have a composite structure to be described in Lemma \ref{lem21}, matching a nonlinear
and essentially parabola-like behavior in a small ball around the origin to an affine linear behavior in a corresponding
outer annulus, the latter increasing so as to coincide with the whole domain $B_R$ at the blow-up time of $\uw$.
The technical challenge, to be addressed in Section \ref{sect3},
will then consist in carefully adjusting the parameters in the definition of $\uw$ in such a manner that
the resulting function in fact has the desired blow-up property, where the cases $n\ge 2$ and $n=1$ will require
partially different arguments (Lemma \ref{lem54} and Lemma \ref{lem55}).
The statement from Theorem \ref{theo14} will thereafter result in Section \ref{sect4}.
\mysection{A parabolic problem satisfied by the mass accumulation function}\label{sect2}
Throughout the sequel, we fix $R>0$ and consider (\ref{0}) in the spatial domain $\Omega:=B_R(0)\subset \R^n$, $n\ge 1$.
Then following a standard procedure (\cite{jaeger_luckhaus}), given a radially symmetric solution $(u,v)=(u(r,t),v(r,t))$
of (\ref{0}) in $\Omega \times [0,T)$ for some $T>0$, we consider the associated mass accumulation function $w$ given by
\be{trans}
	w(s,t):=\int_0^{s^\frac{1}{n}} r^{n-1} u(r,t)dr
	\qquad \mbox{for $s\in [0,R^n]$ and } t\in [0,T).
\ee
In order to describe a basic property of $w$ naturally inherited from $(u,v)$ through (\ref{0}), let us furthermore
introduce the parabolic operator $\parab$ formally given by
\be{parab}
	(\parab \tilde w)(s,t)
	:=\tilde w_t
	- n^2 \cdot \frac{s^{2-\frac{2}{n}} \tilde w_s \tilde w_{ss}}
		{\sqrt{\tilde w_s^2 + n^2 s^{2-\frac{2}{n}} \tilde w_{ss}^2}}
	-n\chi \cdot \frac{(\tilde w-\frac{\mu}{n}s) \cdot \tilde w_s}
		{\sqrt{1+s^{\frac{2}{n}-2} \Big(\tilde w - \frac{\mu}{n}s\Big)^2}}.
\ee
We note here that for $T>0$, the above expression $\parab \tilde w$ is indeed well-defined for all
$t\in (0,T)$ and a.e.~$s\in (0,R^n)$ if, for instance,
$\tilde w \in C^1((0,R^n)\times (0,T))$ is such that $w_s>0$ throughout $(0,R^n)\times (0,T)$ and
$\tilde w(\cdot,t) \in W^{2,\infty}((0,R^n))$ for all $t\in (0,T)$.\\
Now the function $w$ in (\ref{trans}), which clearly complies with these requirements due to smoothness and
positivity of $u$, in fact
solves an appropriate initial-boundary value problem associated with $\parab$:
\begin{lem}\label{lem_w}
  Let $n\ge 1$ and $\chi>0$, and suppose that $(u,v)$ is a positive
  radially symmetric classical
  solution of (\ref{0}) in $\Omega\times (0,T)$ for some $T>0$
  and some nonnegative radially symmetric $u_0\in C^0(\bar\Omega)$.
  Then the function $w$ defined in (\ref{trans}) satisfies
  \be{0w}
	\left\{ \begin{array}{l}
	(\parab w)(s,t)=0, \qquad s\in (0,R^n), \ t\in (0,T), \\[1mm]
	w(0,t)=0, \quad w(R^n,t)=\frac{m}{\omega_n}, \qquad t\in (0,T), \\[1mm]
	w(s,0)=\int_0^{s^\frac{1}{n}} r^{n-1} u_0(r)dr, \qquad s\in [0,R^n],
	\end{array} \right.
  \ee
  where $m:=\io u_0(x)dx$, and where $\omega_n$ denotes the $(n-1)$-dimensional measure of the unit sphere in $\R^n$.
\end{lem}
\proof
  Omitting the arguments $r,t$ and $s:=r^n$ in expressions like $u(r,t)$ and $w(s,t)$,
  upon an integration in the radial version of the first equation in (\ref{0}) we obtain
  \bea{w1}
	w_t
	&=& \int_0^{s^\frac{1}{n}} r^{n-1} u_t(r,t) dr \nn\\
	&=& \int_0^{s^\frac{1}{n}} \bigg\{ \Big(r^{n-1} \frac{uu_r}{\sqrt{u^2+u_r^2}} \Big)_r
	- \chi \Big( r^{n-1} \frac{uv_r}{\sqrt{1+v_r^2}} \Big)_r \bigg\} dr \nn\\
	&=& \big(s^\frac{1}{n}\big)^{n-1} \cdot \frac{uu_r}{\sqrt{u^2+u_r^2}}
	- \chi \cdot \frac{u\cdot (r^{n-1} v_r)}{\sqrt{1+v_r^2}}
  \eea
  for $s\in (0,R^n)$ and $t\in (0,T)$.
  Here in order to replace $v_r$, we integrate the second equation in (\ref{0}), that is, the identity
  $(r^{n-1} v_r)_r=\mu r^{n-1} - r^{n-1}u$, to see that
  \bas
	r^{n-1} v_r=\frac{\mu}{n} \cdot r^n - \int_0^r \rho^{n-1} u(\rho,t) d\rho
	= \frac{\mu}{n}\cdot s - w.
  \eas

Furthermore  (\ref{trans}) can be used to derive
  \bas
	u=nw_s
	\qquad \mbox{and} \qquad
	u_r=n^2 s^{1-\frac{1}{n}} w_{ss},
  \eas
  to infer from (\ref{w1}) that
  \bas
	w_t
	&=& s^{1-\frac{1}{n}} \cdot
	\frac{nw_s \cdot n^2 s^{1-\frac{2}{n}} w_{ss}}{\sqrt{n^2 w_s^2 + n^4 s^{2-\frac{2}{n}} w_{ss}^2}}
	- \chi \cdot \frac{nw_s \cdot (\frac{\mu}{n}s-w)}
		{\sqrt{1+\Big(\frac{\mu}{n} s^\frac{1}{n}-s^{\frac{1}{n}-1} w\Big)^2}} \\
	&=& n^2 \cdot \frac{s^{2-\frac{2}{n}} w_s w_{ss}}{\sqrt{w_s^2 + n^2 s^{2-\frac{2}{n}} w_{ss}^2}}
	+ n\chi \cdot \frac{(w-\frac{\mu}{n}s) \cdot w_s}
		{\sqrt{1+\Big(s^{\frac{1}{n}-1}w - \frac{\mu}{n} s^\frac{1}{n}\Big)^2}}
  \eas
  for $s\in (0,R^n)$ and $t\in (0,T)$.
  This proves the parabolic equation in (\ref{0w}), whereas the statemets therein concerning boundary and initial
  conditions can easily be checked using (\ref{trans}) and the mass conservation property $\io u(x,t)dx=\io u_0(x)dx=m$
  for $t\in (0,T)$.
\mysection{Construction of subsolutions for (\ref{0w})}\label{sect3}
The goal of this section is to construct subsolutions $\uw$ for the parabolic operator introduced in (\ref{parab})
which after some finite time $T$ exhibit a phenomenon of {\em gradient blow-up} in the strong sense that
\bas
	\sup_{s\in (0,R^n)} \frac{\uw(s,t)}{s} \to +\infty
	\qquad \mbox{as } t\nearrow T.
\eas
Since by means of a suitable comparison principle (cf.~Lemma \ref{lem70} in the appendix) we will be able to assert
that $w\ge \uw$ in $[0,R^n] \times [0,T)$, this will entail a similar conclusion for $w$ and hence prove that
$u$ cannot exist as a bounded solution in $\bar\Omega\times [0,T]$.\abs
Our construction will involve several parameters. The first of these is a number $\lambda\in (0,1)$
which eventually, as we shall see later, can be chosen arbitrarily when $n\ge 2$ (see Lemma \ref{lem55}),
but needs to be fixed appropriately close to $1$ in the case $n=1$, depending on the size of the mass $m=\io u_0$
(Lemma \ref{lem54}).
Leaving this final choice open at this point, given any $\lambda\in (0,1)$ we abbreviate
\be{abl}
	\al:=\frac{(1-\lambda)^2}{2\lambda}
	\qquad \mbox{and} \qquad
	\bl:=\frac{3\lambda-1}{2\lambda}
\ee
and introduce
\be{phi}
	\varphi(\xi):=
	\left\{ \begin{array}{ll}
	\lambda \xi^2 \qquad & \mbox{if } \xi \in [0,1], \\[1mm]
	1-\frac{\al}{\xi-\bl} \qquad & \mbox{if } \xi>1.
	\end{array} \right.
\ee
It can then easily be verified that $\varphi$ belongs to $C^1([0,\infty)) \cap W^{2,\infty}((0,\infty)) \cap
C^2([0,\infty) \setminus \{1\})$ with
\be{phi_p}
	\varphi'(\xi)=
	\left\{ \begin{array}{ll}
	2\lambda \xi \qquad & \mbox{if } \xi \in [0,1), \\[1mm]
	\frac{\al}{(\xi-\bl)^2} \qquad & \mbox{if } \xi>1,
	\end{array} \right.
\ee
and
\be{phi_pp}
	\varphi''(\xi)=
	\left\{ \begin{array}{ll}
	2\lambda \qquad & \mbox{if } \xi \in [0,1), \\[1mm]
	- \frac{2 \al}{(\xi-\bl)^3} \qquad & \mbox{if } \xi>1,
	\end{array} \right.
\ee
whence in particular $\varphi'(\xi)>0$ for all $\xi\ge 0$.\abs
With these definitions, we can now specify the basic structure of our comparison functions $\uw$ to be used in the sequel.
Here a second parameter $K$ enters, to be chosen suitably large finally, as well as a parameter function $B$ depending
on time. In combination, these two ingredients determine a line $s=K\sqrt{B(t)}$ in the $(s,t)$-plane which will separate
an inner from an outer region and thereby imply a composite structure of $\uw$ as follows.
\begin{lem}\label{lem21}
  Let $n\ge 1$, $m>0$, $\lambda\in (0,1)$ and $K>1$, and
  suppose that $T>0$ and that $B\in C^1([0,T))$ is such that $B(t)\in (0,1)$ and $K\sqrt{B(t)} < R^n$ for all $t\in [0,T)$
  as well as
  \be{A_well}
	B(t) \le \frac{K^2}{4(\al+\bl)^2}
	\qquad \mbox{for all } t\in [0,T),
  \ee
  where $\al$ and $\bl$ are as in (\ref{abl}).
  Let 	%$\uw:[0,R^n] \times [0,T) \to \R$ be defined by
  \be{uw}
	\uw(s,t):= \left\{ \begin{array}{ll}
	\win(s,t) \qquad & \mbox{if $t\in [0,T)$ and } s\in [0,K\sqrt{B(t)}], \\[1mm]
	\wout(s,t) \qquad & \mbox{if $t\in [0,T)$ and } s\in (K\sqrt{B(t)},R^n],
	\end{array} \right.
  \ee
  where
  \be{win}
	\win(s,t):=
	A(t)\varphi(\xi), \quad \xi=\xi(s,t):=\frac{s}{B(t)},
	\qquad  \mbox{$t\in [0,T)$, } s \in [0, K\sqrt{B(t)}],
  \ee
  with $\varphi$ is as in (\ref{phi}), and where
  \be{wout}
	\wout(s,t):= D(t)s+E(t) \qquad \mbox{for $t\in [0,T)$ and } s\in (K\sqrt{B(t)},R^n],
  \ee
  with
  \be{A}
	A(t)=\frac{m}{\omega_n} \cdot \frac{K^2-2\bl K \sqrt{B(t)} + \bl^2 B(t)}{N(t)}\va \qquad t\in [0,T),
  \ee
%  \be{A_old}
%	A(t)=\frac{m}{\omega_n} \cdot
%	\frac{1-\frac{2\bl}{K} \sqrt{B(t)} + \frac{\bl^2}{K^2} B(t)}
%	{1+\frac{\al R^n}{K^2} - \frac{2(\al+\bl)}{K} \sqrt{B(t)} + \frac{(\al+\bl) \bl}{K^2} B(t)},
%	\qquad t\in [0,T)
%  \ee
  as well as
  \be{D}
	D(t):= \frac{m}{\omega_n} \cdot \frac{\al}{N(t)}\va \qquad t\in [0,T),
  \ee
  and
  \be{E}
	E(t):=\frac{m}{\omega_n} - R^n D(t) \equiv \frac{m}{\omega_n} \cdot \frac{K^2-2(\al+\bl) K\sqrt{B(t)} + (\al+\bl)\bl B(t)}{N(t)}\va
	\qquad t\in [0,T),
  \ee
  with
  \be{N}
	N(t):=K^2 + \al R^n - 2(\al+\bl)K\sqrt{B(t)} + (\al+\bl)\bl B(t),
	\qquad t\in [0,T).
  \ee
  Then $\uw$ is well-defined and continuously differentiable in $[0,R^n]\times [0,T)$ and in addition satisfies
  $\uw(\cdot,t)\in W^{2,\infty}((0,R^n)) \cap C^2([0,R^n] \setminus \{B(t),K\sqrt{B(t)}\})$ for all $t\in [0,T)$
  as well as
  \be{w_bdry}
	w(0,t)=0
	\quad \mbox{and} \quad
	w(R^n,t)=\frac{m}{\omega_n}
	\qquad \mbox{for all } t\in (0,T).
  \ee
  Moreover, the functions $A$ and $D$ defined in (\ref{A}) and (\ref{D}) have their derivatives given by
  \be{A_p}
	A'(t)= \frac{m}{\omega_n} \cdot
	\frac{ \Big(\frac{K}{\sqrt{B(t)}} - \bl\Big) \cdot \Big(\al K^2-\al\bl R^n\Big) \cdot B'(t)}
	{N^2(t)}
  \ee
  and
  \be{D_p}
	D'(t)=\frac{m}{\omega_n} \cdot \frac{\al(\al+\bl) \cdot \Big(\frac{K}{\sqrt{B(t)}}-\bl\Big) \cdot B'(t)}
	{N^2(t)}
  \ee
%  \be{A_p}
%	A'(t)= \frac{m}{\omega_n} \cdot
%	\frac{ \Big(\frac{K}{\sqrt{B(t)}} - \bl\Big) \cdot \Big(\al K^2-\al\bl R^n\Big) \cdot B'(t)}
%	{\Big\{K^2+\al R^n - 2(\al+\bl) K\sqrt{B(t)} +(\al+\bl)\bl B(t) \Big\}^2}
% \ee
%  and
%  \be{D_p}
%	D'(t)=\frac{m}{\omega_n} \cdot \frac{\al(\al+\bl) \cdot \Big(\frac{K}{\sqrt{B(t)}}-\bl\Big) \cdot B'(t)}
%	{\Big\{K^2+\al R^n - 2(\al+\bl) K\sqrt{B(t)} +(\al+\bl)\bl B(t) \Big\}^2}
%  \ee
  for all $t\in (0,T)$.
\end{lem}
\proof
  We first note that for each $t\in [0,T)$, our assumptions that $B(t)>0$ and $K\sqrt{B(t)}<R^n$ ensure
  that both intervals $[0,K\sqrt{B(t)}]$ and $(K\sqrt{B(t)},R^n]$ in (\ref{uw}) are not empty.
  Moreover, thanks to (\ref{A_well}) we have
  \bas
	\frac{2(\al+\bl)}{K}\sqrt{B(t)} \le 1
	\qquad \mbox{for all } t\in [0,T),
  \eas
  which in particular guarantees that the denominators in (\ref{A}), (\ref{D}) and (\ref{E}) are all positive
  and hence $\uw$ well-defined throughout $[0,R^n]\times [0,T)$.
  Moreover, differentiating in (\ref{A}) we can compute
  \bas
	& & \hspace*{-12mm}
	\frac{\omega_n}{m} \cdot N^2(t)A'(t) \\[2mm]
	&=& \bigg\{ -\frac{\bl K B'(t)}{\sqrt{B(t)}} + \bl^2 B'(t) \bigg\} \cdot
	\bigg\{ K^2 + \al R^n - 2(\al+\bl) K\sqrt{B(t)} + (\al+\bl)\bl B(t) \bigg\} \\
	& & - \bigg\{ K^2 + 2\bl K\sqrt{B(t)} + \bl^2 B(t)\bigg\} \cdot
	\bigg\{ -\frac{(\al+\bl) KB'(t)}{\sqrt{B(t)}} + (\al+\bl)\bl B'(t)\bigg\} \\[2mm]
	&=& B'(t) \cdot \bigg\{ \frac{K}{\sqrt{B(t)}} - \bl \bigg\} \cdot
	\Bigg\{ -\bl \cdot \bigg\{ K^2 + \al R^n - 2(\al+\bl)K\sqrt{B(t)} + (\al+\bl)\bl B(t)\bigg\} \\
	& & \hspace*{44mm}
	+ (\al+\bl) \cdot \bigg\{ K^2-2\bl K\sqrt{B(t)} + \bl^2 B(t)\bigg\}
	\Bigg\} \\[2mm]
	&=& B'(t) \cdot \bigg\{ \frac{K}{\sqrt{B(t)}} -\bl\bigg\} \cdot \bigg\{ \al K^2 - \al\bl R^n \bigg\}
	\qquad \mbox{for all } t\in (0,T),
  \eas
  which establishes (\ref{A_p}). Similarly, differentiation in (\ref{D}) readily yields (\ref{D_p}),
  whereas both statements in (\ref{w_bdry}) are direct consequences of (\ref{E}) and the fact that
  $\varphi(0)=0$ according to (\ref{phi}).\abs
  To establish the claimed regularity properties of $\uw$, in view of the above observation that
  $\varphi\in C^1([0,\infty)) \cap W^{2,\infty}((0,\infty)) \cap C^2([0,\infty) \setminus \{1\})$
  we only need to make sure that $\uw$, $\uw_s$ and $\uw_t$ are continuous along the line where $s=K\sqrt{B(t)}$,
  which amounts to showing that
  \be{21.1}
	A(t)\cdot\varphi\Big(\frac{K}{\sqrt{B(t)}}\Big)=D(t)\cdot K\sqrt{B(t)} + E(t)
	\qquad \mbox{for all } t\in [0,T)
  \ee
  and
  \be{21.2}
	\frac{A(t)}{B(t)}\cdot \varphi'\Big(\frac{K}{\sqrt{B(t)}}\Big) = D(t)
	\qquad \mbox{for all } t\in [0,T)
  \ee
  as well as
  \be{21.3}
	A'(t) \cdot \varphi\Big(\frac{K}{\sqrt{B(t)}}\Big)
	- \frac{KA(t)B'(t)}{\sqrt{B(t)}^3} \cdot \varphi' \Big(\frac{K}{\sqrt{B(t)}}\Big)
	= D'(t) \cdot K\sqrt{B(t)} + E'(t)
	\qquad \mbox{for all } t\in [0,T).
  \ee
  To derive (\ref{21.1}), we use (\ref{phi}) to see that
  \bas
	& & \hspace*{-20mm}
	A(t)\cdot\varphi\Big(\frac{K}{\sqrt{B(t)}}\Big) + \Big(R^n-K\sqrt{B(t)}\Big) \cdot  D(t) \\[1mm]
	&=& \frac{m}{\omega_n} \cdot
	\frac{(K-\bl \sqrt{B(t)})^2}
	{K^2+\al R^n - 2(\al+\bl) K\sqrt{B(t)} + (\al+\bl)\bl B(t)}
	\cdot \frac{\frac{K}{\sqrt{B(t)}}-\al-\bl}{\frac{K}{\sqrt{B(t)}}-\bl} \\
	& & + (R^n-K\sqrt{B(t)}) \cdot \frac{m}{\omega_n} \cdot
	\frac{a}
	{K^2+\al R^n - 2(\al+\bl) K\sqrt{B(t)} + (\al+\bl)\bl B(t)} \\[2mm]
	&=& \frac{m}{\omega_n} \cdot
	\frac{ \frac{K-(\al+\bl) \sqrt{B(t)}}{K-\bl \sqrt{B(t)}} \cdot (K-\bl \sqrt{B(t)})^2
	+ \al (R^n-K\sqrt{B(t)})}
	{K^2+\al R^n - 2(\al+\bl) K\sqrt{B(t)} + (\al+\bl)\bl B(t)} \\[2mm]
	&=& \frac{m}{\omega_n} \cdot
	\frac{K^2 - \bl K\sqrt{B(t)} -(\al+\bl)K\sqrt{B(t)} + (\al+\bl)\bl B(t) + \al R^n - \al K\sqrt{B(t)}}
	{K^2+\al R^n - 2(\al+\bl) K\sqrt{B(t)} + (\al+\bl)\bl B(t)} \\[2mm]
	&=& \frac{m}{\omega_n}
	\qquad \mbox{for all } t\in [0,T),
  \eas
  which due to (\ref{E}) means that indeed
  \bas
	A(t)\cdot\varphi\Big(\frac{K}{\sqrt{B(t)}}\Big)-D(t)\cdot K\sqrt{B(t)} - E(t)
	= A(t)\cdot\varphi\Big(\frac{K}{\sqrt{B(t)}}\Big) + (R^n-K\sqrt{B(t)}) D(t) - \frac{m}{\omega_n}
	= 0
  \eas
  for all $t\in [0,T)$.
  Next, from (\ref{D}) and (\ref{A}) it immediately follows that
  \bas
	\frac{D(t)}{A(t)}
	&=& \frac{1}{K^2} \cdot \frac{\al}{1-\frac{2\bl}{K}\sqrt{B(t)} + \frac{\bl^2}{K^2} B(t)} \\
	&=& \frac{\al}{K^2 - 2\bl K \sqrt{B(t)} + \bl^2 B(t)} \\
	&=& \frac{\al}{B(t)\Big(\frac{K}{\sqrt{B(t)}} - \bl\Big)^2} \\
	&=& \frac{1}{B(t)} \cdot \varphi'\Big(\frac{K}{\sqrt{B(t)}}\Big)
	\qquad \mbox{for all } t\in [0,T),
  \eas
  which establishes (\ref{21.2}).\\
  Finally, in verifying (\ref{21.3}) we make use of (\ref{21.1}) and (\ref{21.2}) as well as
  (\ref{A}), (\ref{D}), (\ref{E}), (\ref{A_p}) and (\ref{D_p})
  to see that
  \bas
	& & \hspace*{-10mm}
	A'(t) \cdot \varphi\Big(\frac{K}{\sqrt{B(t)}}\Big)
	- \frac{KA(t)B'(t)}{\sqrt{B(t)}^3} \cdot \varphi' \Big(\frac{K}{\sqrt{B(t)}}\Big)
	- D'(t) \cdot K\sqrt{B(t)} + E'(t) \\[2mm]
	&=& A'(t) \cdot \frac{\frac{m}{\omega_n} - \Big(R^n-K\sqrt{B(t)}\Big)\cdot D(t)}{A(t)}
	- \frac{K A(t) B'(t)}{\sqrt{B(t)}^3} \cdot \frac{B(t) D(t)}{A(t)}
	+ \Big(R^n - K\sqrt{B(t)} \Big)\cdot D'(t) \\[2mm]
	&=& \frac{m}{\omega_n} \cdot
	\frac{\Big(\frac{K}{\sqrt{B(t)}} - \bl\Big) \cdot (\al K^2 - \al \bl R^n) \cdot B'(t)}{N^2(t)}
	\cdot \bigg\{ \frac{m}{\omega_n} \cdot \frac{K^2 -2\bl K\sqrt{B(t)} + \bl^2 B(t)}{N(t)} \bigg\}^{-1}
	\times \\
	& & \hspace*{10mm}
	\times \bigg\{
	\frac{m}{\omega_n} - \Big(R^n-K\sqrt{B(t)}\Big) \cdot \frac{m}{\omega_n} \cdot \frac{\al}{N(t)}
	\bigg\} \\[1mm]
	& & - \frac{KB'(t)}{\sqrt{B(t)}} \cdot \frac{m}{\omega_n} \cdot \frac{\al}{N(t)} \\
	& & + \Big(R^n-K\sqrt{B(t)}\Big) \cdot \frac{m}{\omega_n} \cdot
	\frac{\al(\al+\bl) \Big(\frac{K}{\sqrt{B(t)}}-\bl\Big) \cdot B'(t)}{N^2(t)} \\[2mm]
	&=& \frac{m}{\omega_n} \cdot \frac{\al B'(t)}{\sqrt{B(t)} N^2(t)} \cdot
	\Bigg\{ \frac{ \Big(K-\bl\sqrt{B(t)} \Big) \cdot \Big(K^2-\bl R^n\Big)}{\Big(K-b\sqrt{B(t)}\Big)^2}
	\cdot \bigg[ N(t)-\Big(R^n-K\sqrt{B(t)}\Big) \cdot \al\bigg] \\
	& & \hspace*{35mm}
	-KN(t)
	+ (\al+\bl) \Big(R^n-K\sqrt{B(t)}\Big) \cdot \Big(K-\bl \sqrt{B(t)}\Big)
	\Bigg\} \\[2mm]
	&=& \frac{m}{\omega_n} \cdot
	\frac{\al B'(t)}{\Big(K-\bl \sqrt{B(t)}\Big)\sqrt{B(t)} N^2(t)} \cdot
	\Bigg\{ \Big(K^2-\bl R^n\Big) \cdot \Big[ K^2-(\al+2\bl) K\sqrt{B(t)} + (\al+\bl)\bl B(t)\Big] \\
	& & \hspace*{40mm}
	- K\Big(K-b\sqrt{B(t)}\Big) \cdot \Big[ K^2+ \al R^n - 2(\al+\bl) K\sqrt{B(t)} + (\al+\bl)\bl B(t)\Big] \\
	& & \hspace*{40mm}
	+ (\al+\bl) \Big(R^n-K\sqrt{B(t)}\Big) \cdot \Big(K-\bl \sqrt{B(t)}\Big)^2 \Bigg\}
  \eas
  for all $t\in (0,T)$.
  Since it can be checked in a straightforward manner that herein we have
  \bas
	& & \hspace*{-20mm}
	\Big(K^2-\bl R^n\Big) \cdot \Big[ K^2-(\al+2\bl) K\sqrt{B(t)} + (\al+\bl)\bl B(t)\Big] \\
	& &
	- K\Big(K-b\sqrt{B(t)}\Big) \cdot \Big[ K^2+ \al R^n - 2(\al+\bl) K\sqrt{B(t)} + (\al+\bl)\bl B(t)\Big] \\
	& &
	+ (\al+\bl) \Big(R^n-K\sqrt{B(t)}\Big) \cdot \Big(K-\bl \sqrt{B(t)}\Big)^2 \\[2mm]
	&=& 0,
  \eas
  this shows (\ref{21.3}) and thereby completes the proof.
\qed
\subsection{Subsolution properties: Outer region}
Let us first make sure that if the function $B$ entering the above definition of $\uw$ is suitably small and satisfies
an appropriate differential inequality, then $\uw$ becomes a subsolution in the corresponding outer region addressed
in (\ref{uw}).
\begin{lem}\label{lem60}
  Let $n\ge 1, \chi>0, m>0, \lambda \in (0,1), K>1$ and $B_0\in (0,1)$ be such that $K\sqrt{B_0}<R^n$ and
  \be{60.1}
	B_0 \le \frac{K^2}{16(\al+\bl)^2}
  \ee
  with $\al$ and $\bl$ given by (\ref{abl}).
  Then if for some $T>0$, $B\in C^1([0,T))$ is positive and nonincreasing and such that
  \be{60.2}
	\left\{ \begin{array}{l}
	B'(t) \ge - \frac{nm\chi K}{2(\al+\bl) \omega_n R^n \sqrt{1+K^{\frac{2}{n}-2} \frac{m^2}{\omega_n^2}}} \cdot
	B^{1-\frac{1}{2n}}(t),
	\qquad t\in (0,T), \\[1mm]
	B(0) \le B_0,
	\end{array} \right.
  \ee
  the function $\wout$ defined in (\ref{wout}) satisfies
  \be{60.3}
	(\parab \wout)(s,t) \le 0
	\qquad \mbox{for all $t\in (0,T)$ and } s\in (K\sqrt{B(t)}, R^n)
  \ee
  with $\parab$ given by (\ref{parab}).
\end{lem}
\proof
  Again using that $E(t) =\frac{m}{\omega_n} - R^n D(t)$ for all $t\in (0,T)$ by (\ref{E}), we have
  \be{60.4}
	\wout(s,t) = D(t)s + E(t)
	=\frac{m}{\omega_n} - D(t) \cdot (R^n-s)
	\qquad \mbox{for all $t\in (0,T)$ and } s\in (K\sqrt{B(t)}, R^n),
  \ee
  so that recalling (\ref{D_p}) we obtain
  \bea{60.5}
	(\wout)_t(s,t)
	&=& - D'(t) \cdot (R^n-s) \nn\\
	&=& - \frac{m}{\omega_n} \cdot
	\frac{\al (\al+\bl) \cdot \Big\{ \frac{K}{\sqrt{B(t)}} - \bl \Big\}}{N^2(t)}
%	{\Big\{K^2 + \al R^n - 2(\al+\bl) K\sqrt{B(t)} + (\al+\bl)\bl B(t) \Big\}^2}
	\cdot B'(t) \cdot (R^n-s)
%	\qquad \mbox{for all $t\in (0,T)$ and } s\in (K\sqrt{B(t)}, R^n)
  \eea
  for all $t\in (0,T)$ and $s\in (K\sqrt{B(t)},R^n)$, where $N(t)$ is as in (\ref{N}) for such $t$.
%abbreviated
%  \be{60.555}
%	N(s,t):=\Big\{K^2 + \al R^n - 2(\al+\bl) K\sqrt{B(t)} + (\al+\bl)\bl B(t) \Big\}^2
%	\qquad t\in (0,T), \ s\in (K\sqrt{B(t)},R^n).
%  \ee
%  for such $t$ and $s$.
%  \bea{60.5}
%	& & \hspace*{-20mm}
%	(\wout)_t(s,t)
%	\ = \ - D'(t) \cdot (R^n-s) \nn\\
%	& & \hspace*{-10mm}
%	= \ - \frac{m}{\omega_n} \cdot
%	\frac{\al \cdot \Big\{ \frac{(\al+\bl)K}{\sqrt{B(t)}} - (\al+\bl)\bl \Big\}}
%	{\Big\{K^2 + \al R^n - 2(\al+\bl) K\sqrt{B(t)} + (\al+\bl)\bl B(t) \Big\}^2}
%	\cdot B'(t) \cdot (R^n-s)
%  \eea
%  for all $t\in (0,T)$ and $s\in (K\sqrt{B(t)},R^n)$.
  In order to compensate the positive contribution of this term $(\wout)_t$ to $\parab \wout$ by a
  suitably negative impact of the rightmost term
  \be{60.55}
	I(s,t)
	:= -n\chi \cdot \frac{(\wout-\frac{\mu}{n}s) \cdot (\wout)_s}{\sqrt{1+s^{\frac{2}{n}-2} (\wout-\frac{\mu}{n}s)^2}},
	\qquad t\in (0,T), \ s\in (K\sqrt{B(t)}, R^n),
  \ee
  in (\ref{parab}), we use (\ref{60.4}) and (\ref{D}) to rewrite
  \bea{60.6}
	\wout(s,t)-\frac{\mu}{n}s
	&=& \Big(\frac{m}{\omega_n R^n} - D(t)\Big) \cdot (R^n-s) \nn\\[1mm]
	&=& \frac{m}{\omega_n R^n} \cdot \Big(1-R^n \cdot \frac{\omega_n}{m} \cdot D(t)\Big) \cdot (R^n-s) \nn\\[1mm]
	&=& \frac{m}{\omega_n R^n} \cdot
	\bigg( 1-\frac{\al R^n}
	{K^2+ \al R^n - 2(\al+\bl) K\sqrt{B(t)} + (\al+\bl)\bl B(t)}
	\bigg) \cdot (R^n-s) \nn\\[1mm]
	&=& \frac{m}{\omega_n R^n} \cdot
	\frac{K^2-2(\al+\bl) K\sqrt{B(t)} + (\al+\bl)\bl B(t)}
	{K^2+ \al R^n - 2(\al+\bl) K\sqrt{B(t)} + (\al+\bl)\bl B(t)}
	\cdot (R^n-s)
  \eea
  for $t\in (0,T)$ and $s\in (K\sqrt{B(t)}, R^n)$.
  As
  \be{60.66}
	K^2-2(\al+\bl) K\sqrt{B(t)} \ge \frac{1}{2} K^2
	\qquad \mbox{for all } t\in (0,T)
  \ee
  by (\ref{60.1}), this in particular implies that
  \bas
	0 < \wout(s,t) - \frac{\mu}{n}s \le \frac{m}{\omega_n R^n} \cdot R^n = \frac{m}{\omega_n}
  \eas
  and hence
  \bas
	\Big(\wout(s,t)-\frac{\mu}{n}s\Big)^2 \le \frac{m^2}{\omega_n^2}
	\qquad \mbox{for all $t\in (0,T)$ and } s\in (K\sqrt{B(t)}, R^n).
  \eas
  Since moreover $1\le B^{\frac{1}{n}-1}(t)$ for all $t\in (0,T)$ due to the fact that $B_0<1$, we can thus estimate
  the denominator in (\ref{60.55}) in the considered outer region according to
  \bea{60.7}
	\sqrt{1-s^{\frac{2}{n}-2} \Big(\wout(s,t)-\frac{\mu}{n}s\Big)^2}
	&\le& \sqrt{ B^{\frac{1}{n}-1}(t) + \Big(K\sqrt{B(t)}\Big)^{\frac{2}{n}-2} \cdot \frac{m^2}{\omega_n^2}} \nn\\
	&=& \sqrt{1+K^{\frac{2}{n}-2} \cdot \frac{m^2}{\omega_n^2}} \cdot B^{\frac{1}{2n}-\frac{1}{2}}(t) \nn\\[2mm]
	& & \hspace*{10mm} \qquad \mbox{for all $t\in (0,T)$ and } s\in (K\sqrt{B(t)}, R^n).
  \eea
  Using (\ref{60.7}) and (\ref{60.6}) and that
  \bas
	(\wout)_s(s,t)
	= D(t)
	= \frac{m}{\omega_n} \cdot \frac{\al}{K^2+\al R^n - 2(\al+\bl) K\sqrt{B(t)} + (\al+\bl)\bl B(t)}
%	\qquad \mbox{for all $t\in (0,T)$ and } s\in (K\sqrt{B(t)}, R^n),
  \eas
  for $t\in (0,T)$ and $s\in (K\sqrt{B(t)}, R^n)$,
  we thereby find that
  \bas
	-I(s,t)
	&\ge& n\chi \cdot \frac{\Big(\wout(s,t)-\frac{\mu}{n}s\Big) \cdot (\wout)_s(s,t)}
	{\sqrt{1+K^{\frac{2}{n}-2} \cdot \frac{m^2}{\omega_n^2}} \cdot B^{\frac{1}{2n}-\frac{1}{2}}(t)} \\[2mm]
	&=& \frac{n\chi}{\sqrt{1+K^{\frac{2}{n}-2} \cdot \frac{m^2}{\omega_n^2}}} \cdot B^{\frac{1}{2}-\frac{1}{2n}}(t)
	\times \\
	& & \hspace*{20mm}
	\times \frac{m}{\omega_n R^n} \cdot
	\frac{K^2-2(\al+\bl) K\sqrt{B(t)} + (\al+\bl)\bl B(t)}
	{K^2+ \al R^n - 2(\al+\bl) K\sqrt{B(t)} + (\al+\bl)\bl B(t)}
	\cdot (R^n-s) \times \\
	& & \hspace*{20mm}
	\times \frac{m}{\omega_n} \cdot \frac{\al}{K^2+\al R^n - 2(\al+\bl) K\sqrt{B(t)} + (\al+\bl)\bl B(t)} \\[2mm]
	&=&  \frac{m}{\omega_n} \cdot \frac{\al}{N^2(t)}
	\cdot (R^n-s) \cdot c_1 \cdot \Big\{ K^2 - 2(\al+\bl) K\sqrt{B(t)} + (\al+\bl)\bl B(t)\Big\} \cdot
	B^{\frac{1}{2}-\frac{1}{2n}}(t)
%	&=&  \frac{m}{\omega_n} \cdot \frac{\al}{\Big\{ K^2+ \al R^n - 2(\al+\bl) K\sqrt{B(t)} + (\al+\bl)\bl B(t) \Big\}^2}
%	\cdot (R^n-s) \cdot c_1 \cdot \Big\{ K^2 - 2(\al+\bl) K\sqrt{B(t)} + (\al+\bl)\bl B(t)\Big\} \cdot
%	B^{\frac{1}{2}-\frac{1}{2n}}(t)
%	\qquad \mbox{for all $t\in (0,T)$ and } s\in (K\sqrt{B(t)}, R^n)
  \eas
  for all $t\in (0,T)$ and $s\in (K\sqrt{B(t)}, R^n)$,
  with $N$ as defined in (\ref{N}) and
  \bas
	c_1:=\frac{nm\chi}{\omega_n R^n \cdot \sqrt{1+K^{\frac{2}{n}-2} \cdot \frac{m^2}{\omega_n^2}}}.
  \eas
  Since evidently $(\wout)_{ss} \equiv 0$, combining this with (\ref{60.66}) and (\ref{60.5}) shows that
  \bea{60.8}
	(\parab \wout)(s,t)
	&\le&  \frac{m}{\omega_n} \cdot \frac{\al}{N^2(t)}
	%\Big\{ K^2+ \al R^n - 2(\al+\bl) K\sqrt{B(t)} + (\al+\bl)\bl B(t) \Big\}^2}
	\cdot (R^n-s) \cdot \Bigg\{
	\Big\{ -\frac{(\al+\bl)K}{\sqrt{B(t)}} + (\al+\bl)\bl \Big\} \cdot B'(t) \nn\\
	& & \hspace*{30mm}
	- c_1 \cdot \Big\{ K^2-2(\al+\bl) K\sqrt{B(t)} + (\al+\bl)\bl B(t)\Big\} \cdot B^{\frac{1}{2}-\frac{1}{2n}}(t) \Bigg\}
	\nn\\[2mm]
	&\le&  \frac{m}{\omega_n} \cdot \frac{\al}{N^2(t)}
	%\Big\{ K^2+ \al R^n - 2(\al+\bl) K\sqrt{B(t)} + (\al+\bl)\bl B(t) \Big\}^2}
	\cdot (R^n-s) \cdot \bigg\{
	-\frac{(\al+\bl)K}{\sqrt{B(t)}} \cdot B'(t)
	-\frac{c_1}{2} K^2 B^{\frac{1}{2}-\frac{1}{2n}}(t)\bigg\}
%	\qquad \mbox{for all $t\in (0,T)$ and } s\in (K\sqrt{B(t)}, R^n),
  \eea
  for all $t\in (0,T)$ and $s\in (K\sqrt{B(t)}, R^n)$,
  because $(\al+\bl)\bl B'(t) \le 0$ for all $t\in (0,T)$.
  In view of the definition of $c_1$, (\ref{60.2}) warrants that herein
  \bas
	-\frac{(\al+\bl)K}{\sqrt{B(t)}} \cdot B'(t)
	-\frac{c_1}{2} K^2 B^{\frac{1}{2}-\frac{1}{2n}}(t)
	&=& \frac{(\al+\bl)K}{\sqrt{B(t)}} \cdot \bigg\{ -B'(t) - \frac{c_1 K}{2(\al+\bl)} \cdot B^{1-\frac{1}{2n}}(t)\bigg\}
	\\[2mm]
	&\le& 0
	\qquad \mbox{for all $t\in (0,T)$ and } s\in (K\sqrt{B(t)}, R^n),
  \eas
  so that (\ref{60.3}) results from (\ref{60.8}).
\qed
\subsection{Subsolution properties: Inner region}
We proceed to study under which assumptions on the parameters the function $\uw$ defines a subsolution
in the corresponding inner domain.
To prepare our analysis, let us first compute the action of the operator $\parab$ on $\uw$ in the respective
region as follows.
\begin{lem}\label{lem235}
  Let $n\ge 1, \chi>0$, $m>0$, $\lambda\in (0,1), K>1$ and $T>0$, and suppose that
  $B\in C^1([0,T))$ is positive and satisfies (\ref{A_well}) as well as $K\sqrt{B(t)}<R^n$ for all $t\in [0,T)$.
  Then the function $\win$ defined in (\ref{win}) has the property that
  \bea{pin}
%	\frac{B(t)}{A(t)\varphi'(\xi)} \cdot
	(\parab \win)(s,t)
	&=& A'(t)\varphi(\xi)
	+ \frac{A(t)\varphi'(\xi)}{B(t)} \cdot \Big\{
	-\xi B'(t) + J_1(s,t) + J_2(s,t) \Big\} \nn\\[2mm]
	& & \hspace*{20mm}
	\qquad \mbox{for all $t\in (0,T)$ and } s\in (0,K\sqrt{B(t)}) \setminus \{B(t)\},
  \eea
  where $\xi=\xi(s,t)=\frac{s}{B(t)}$, $\parab$ is as in (\ref{parab}) and
  \be{j1}
	J_1(s,t):=-n^2 \cdot \frac{\xi^{2-\frac{2}{n}}\varphi''(\xi)}
	{\sqrt{B^{\frac{4}{n}-2}(t) \varphi'^2(\xi) + n^2 B^{\frac{2}{n}-2}(t) \xi^{2-\frac{2}{n}} \varphi''^2(\xi)}}
  \ee
  and
  \be{j2}
	J_2(s,t):=-n\chi \cdot \frac{A(t)\varphi(\xi)-\frac{\mu}{n} B(t)\xi}
	{\sqrt{1+B^{\frac{2}{n}-2}(t)\xi^{\frac{2}{n}-2} \cdot \Big(A(t)\varphi(\xi)-\frac{\mu}{n}B(t)\xi\Big)^2}}
  \ee
%  with $\xi=\xi(s,t)=\frac{s}{B(t)}$
  for $t\in (0,T)$ and $s\in (0,K\sqrt{B(t)}) \setminus \{B(t)\}$.
\end{lem}
\proof
  Since $\xi_t=-\frac{sB'(t)}{B^2(t)}=-\frac{\xi B'(t)}{B(t)}$ and $\xi_s=\frac{1}{B(t)}$,
  %suppressing the argument $\xi$ in $\varphi(\xi)$ etc.~
  we can compute
  \be{235.1}
	(\win)_t=A'(t)\varphi(\xi) - \frac{A(t)\xi B'(t)}{B(t)} \cdot \varphi'(\xi)
  \ee
  as well as
  \bas
	(\win)_s=\frac{A(t)}{B(t)} \cdot \varphi'(\xi)
	\quad \mbox{and} \quad
	(\win)_{ss}=\frac{A(t)}{B^2(t)} \cdot \varphi''(\xi)
  \eas
  for all $t\in (0,T)$ and $s\in (0,K\sqrt{B(t)}) \setminus \{B(t)\}$.
  Therefore,
  \bea{235.2}
	n^2 \cdot \frac{s^{2-\frac{2}{n}} (\win)_s (\win)_{ss}}{\sqrt{(\win)_s^2 + n^2 s^{2-\frac{2}{n}} (\win)_{ss}^2}}
	&=& n^2 \cdot
	\frac{(B(t)\xi)^{2-\frac{2}{n}} \cdot \frac{A(t)}{B(t)} \varphi'(\xi) \cdot \frac{A(t)}{B^2(t)} \varphi''(\xi)}
	{\sqrt{\Big(\frac{A(t)}{B(t)}\varphi'(\xi)\Big)^2 + n^2 (B(t)\xi)^{2-\frac{2}{n}} \cdot
		\Big(\frac{A(t)}{B^2(t)} \varphi''(\xi)\Big)^2}} \nn\\
	&=& \frac{n^2 A(t)\varphi'(\xi)}{B(t)} \cdot
	\frac{\xi^{2-\frac{2}{n}} \varphi''(\xi)}
	{\sqrt{B^{\frac{4}{n}-2}(t) \varphi'^2(\xi) + n^2 B^{\frac{2}{n}-2}(t)\xi^{2-\frac{2}{n}} \varphi''^2(\xi)}} \nn\\
	&=& \frac{A(t)\varphi'(\xi)}{B(t)} \cdot J_1(s,t)
  \eea
  and
  \bea{235.3}
	n\chi \cdot \frac{ \Big(\win-\frac{\mu}{n}s\Big) (\win)_s}
	{\sqrt{1+s^{\frac{2}{n}-2} \Big(\win-\frac{\mu}{n}s\Big)^2}}
	&=& n\chi \cdot \frac{\Big(A(t)\varphi(\xi)-\frac{\mu}{n} \cdot B(t)\xi \Big)\cdot \frac{A(t)}{B(t)}\varphi'(\xi)}
	{\sqrt{1+(B(t)\xi)^{\frac{2}{n}-2} \cdot \Big(A(t)\varphi(\xi)-\frac{\mu}{n} \cdot B(t)\xi\Big)^2}} \nn\\[2mm]
	&=&\frac{A(t)\varphi'(\xi)}{B(t)} \cdot J_2(s,t)
  \eea
  for any such $t$ and $s$.
  By definition (\ref{parab}) of $\parab$, (\ref{235.1})-(\ref{235.3}) prove (\ref{pin}).
\qed
In further examining (\ref{pin}), it will be convenient to know that the factor $A$ appearing in (\ref{win}) is
nonincreasing with time, meaning that the first summand on the right-hand side in (\ref{pin}) will be nonpositive.
It is the objective of the following lemma to assert that this can indeed be achieved by choosing
the function $B$ to be nonincreasing and appropriately small throughtout $[0,T)$.
\begin{lem}\label{lem44}
  Let $n\ge 1, m>0, \lambda\in (0,1)$ and $K>1$ be such that $K\ge\sqrt{\bl R^n}$, and suppose that $B_0\in (0,1)$
  satisfies
  \be{44.01}
	B_0 \le \frac{K^2}{4(\al+\bl)^2}.
%	B_0\le \frac{K^2}{\bl^2}.
  \ee
  Then if $T>0$ and $B\in C^1([0,T))$ is a positive and nonincreasing function fulfilling $B(0)\le B_0$,
  for the function $A$ in (\ref{A}) we have
  \be{44.1}
	A'(t) \le 0
	\qquad \mbox{for all } t\in (0,T).
  \ee
  In particular,
  \be{AT}
	A(t) \ge A_T:=\frac{m}{\omega_n} \cdot \frac{1}{1+\frac{\al R^n}{K^2}}
	\qquad \mbox{for all } t\in (0,T).
  \ee
\end{lem}
\proof
  We recall that by (\ref{A_p}), with $N$ given by (\ref{N}) we have
  \be{44.3}
	\frac{\omega_n}{m} \cdot N^2(t)A'(t)
	= \Big( \frac{K}{\sqrt{B(t)}} -\bl\Big)  \cdot \Big( \al K^2 - \al\bl R^n \Big) \cdot B'(t)
	\qquad \mbox{for all } t\in (0,T).
  \ee
  Here since our assumption (\ref{44.01}) implies that $B_0 \le \frac{K^2}{\bl^2}$, by monotonicity of $B$ we obtain that
  \bas
	\frac{K}{\sqrt{B(t)}} - \bl \ge \frac{K}{\sqrt{B_0}} - \bl \ge 0
	\qquad \mbox{for all } t\in (0,T),
  \eas
  whereas the inequality $K\ge \sqrt{\bl R^n}$ ensures that
  \bas
	\al K^2 - \al \bl R^n \ge 0.
  \eas
  Again using that $B' \le 0$, from (\ref{44.3}) we thus conclude that (\ref{44.1}) holds, whereupon (\ref{AT}) follows
  upon taking $t\nearrow T$ in (\ref{A}).
\qed
\subsection{Subsolution properties: Very inner region}
Now in the part very near the origin where $s<B(t)$ and hence $\xi=\frac{s}{B(t)}<1$, the expression
$J_2$ in (\ref{pin}), originating from the chemotactic term in (\ref{0}), need not be positive due to (\ref{phi})
and the linear growth of the minuend $\frac{\mu}{n}B(t)\xi$ in the numernator in (\ref{j2}).
Fortunately, it turns out that the respective unfavorable effect of this to
$\parab \win$ in (\ref{pin}) can be overbalanced by a suitable contribution of $J_1$,
which in fact is negative in this region due to the convexity of $\varphi$ on $(0,1)$.
Under an additional smallness assumption on $B$, we can indeed achieve the following.
\begin{lem}\label{lem61}
  Let $n\ge 1, \chi>0, m>0, \lambda\in (0,1), K>1$ and $B_0\in (0,1)$ be such that $K\sqrt{B_0}<R^n$ and
  \be{61.111}
	B_0 \le \frac{K^2}{4(\al+\bl)^2}
  \ee
  as well as
  \be{61.1}
	B_0 \le \Big(\frac{n}{4\chi\mu}\Big)^n.
  \ee
  Suppose that $T>0$, and that $B\in C^1([0,T))$ is a positive and nonincreasing function satisfying
  \be{61.2}
	\left\{ \begin{array}{l}
	B'(t) \ge - \frac{n}{4} B^{1-\frac{1}{n}}(t), \qquad t\in (0,T), \\[1mm]
	B(0) \le B_0.
	\end{array} \right.
  \ee
  Then the function $\win$ defined in (\ref{win}) has the property that
  \be{61.3}
	(\parab \win)(s,t) \le 0
	\qquad \mbox{for all $t\in (0,T)$ and } s\in (0,B(t)).
  \ee
\end{lem}
\proof
  Writing $\xi=\frac{s}{B(t)}$ for $t\in (0,T)$ and $s\in (0,B(t))$,
  in (\ref{pin}) we can estimate the taxis term from above according to
  \bea{61.4}
	J_2(s,t)
	&=& - n\chi \cdot \frac{A(t)\varphi(\xi)-\frac{\mu}{n}B(t)\xi}
	{\sqrt{1+B^{\frac{2}{n}-2}(t) \xi^{\frac{2}{n}-2} \Big(A(t)\varphi(\xi)-\frac{\mu}{n}B(t)\xi \Big)^2}} \nn\\
	&\le& n\chi \cdot \frac{\frac{\mu}{n}B(t)\xi}{\sqrt{1}} \nn\\[2mm]
	&=& \chi\mu B(t)\xi
	\qquad \mbox{for all $t\in (0,T)$ and } s\in (0,B(t)).
  \eea
%  The effect on this to $\parab \win$ in (\ref{pin}) can be overbalanced by using the convexity of $\varphi$
%  on $(0,1)$. To see this, let us
  We next recall that since $\xi\in (0,1)$ whenever $s\in (0,B(t))$, and hence
  $\varphi'(\xi)=2\lambda\xi$ and $\varphi''(\xi)=2\lambda$, we have
  \bas
	\frac{B^{\frac{4}{n}-2}(t) \varphi'^2(\xi)}{B^{\frac{2}{n}-2} \xi^{\frac{2}{n}-2} \varphi''^2(\xi)}
	&=& B^\frac{2}{n}(t) \cdot \frac{\varphi'^2(\xi)}{\xi^{\frac{2}{n}-2} \varphi''^2(\xi)} \\
	&=& B^\frac{2}{n}(t) \cdot \frac{4\lambda^2 \xi^2}{\xi^{\frac{2}{n}-2} \cdot 4\lambda^2} \\
	&=& B^\frac{2}{n}(t) \cdot \xi^{4-\frac{2}{n}} \\
	&\le& B^\frac{2}{n}(t) \\
	&\le& 1
	\qquad \mbox{for all $t\in (0,T)$ and } s\in (0,B(t)),
  \eas
  because $B \le B_0 \le 1$ throughout $(0,T)$.
  Now since $\sqrt{\varphi''^2(\xi)}=\varphi''(\xi)$ thanks to the convexity of $\varphi$ on $(0,1)$,
  in (\ref{j1}) we therefore find that
  \bas
	-J_1(s,t)
	&=& n^2 \cdot \frac{\xi^{2-\frac{2}{n}}\varphi''(\xi)}
	{\sqrt{B^{\frac{4}{n}-2}(t) \varphi'^2(\xi) + n^2 B^{\frac{2}{n}-2}(t) \xi^{\frac{2}{n}-2}
		\varphi''^2(\xi)}} \\
	&\ge& n^2 \cdot \frac{\xi^{2-\frac{2}{n}} \varphi''(\xi)}
	{\sqrt{(1+n^2) B^{\frac{2}{n}-2}(t) \xi^{\frac{2}{n}-2} \varphi''^2(\xi)}} \\
	&=& \frac{n^2}{\sqrt{1+n^2}} \cdot B^{1-\frac{1}{n}}(t) \xi^{1-\frac{1}{n}}
	\qquad \mbox{for all $t\in (0,T)$ and } s\in (0,B(t)).
  \eas
  As $\sqrt{1+n^2} \le 2n$ and hence $\frac{n^2}{\sqrt{1+n^2}} \ge \frac{n}{2}$, due to (\ref{61.4})
  we thereby obtain from (\ref{pin}), applying Lemma \ref{lem44} on the basis of (\ref{61.111}), that
  \bea{61.5}
	\frac{B(t)}{A(t)\varphi'(\xi)} \cdot (\parab \win)(s,t)
	&\le& -\xi B'(t) - \frac{n}{2} B^{1-\frac{1}{n}}(t) \xi^{1-\frac{1}{n}} + \chi\mu B(t)\xi \nn\\
	&=& - \xi \cdot \Big\{ -B'(t) -\frac{n}{2} B^{1-\frac{1}{n}}(t) \xi^{-\frac{1}{n}} + \chi\mu B(t) \Big\}
%	\qquad \mbox{for all $t\in (0,T)$ and } s\in (0,B(t)).
  \eea
  for all $t\in (0,T)$ and $s\in (0,B(t))$.
  Here, using that $\xi<1$ implies that $\xi^{-\frac{1}{n}}\ge 1$, and that the restriction (\ref{61.1}) on
  $B_0$ ensures that
  \bas
	\frac{\chi\mu B(t)}{\frac{n}{4} B^{1-\frac{1}{n}}(t)}
	= \frac{4\chi\mu}{n} \cdot B^\frac{1}{n}(t)
	\le 1,
  \eas
  we see that
  \bas
	-B'(t) - \frac{n}{2} B^{1-\frac{1}{n}}(t) \xi^{-\frac{1}{n}} + \chi\mu B(t)
	&\le& - B'(t) - \frac{n}{2} B^{1-\frac{1}{n}}(t) + \frac{n}{4} B^{1-\frac{1}{n}}(t) \\
	&=& -B'(t) - \frac{n}{4} B^{1-\frac{1}{n}}(t)
	\qquad \mbox{for all $t\in (0,T)$ and } s\in (0,B(t)).
  \eas
  As a consequence of (\ref{61.2}), the claim therefore results from (\ref{61.5}).
\qed
\subsection{Subsolution properties: Intermediate region}
The crucial part of our analysis will be concerned with the remaining intermediate region, that is, the
outer part of the inner domain where $B(t)<s<K\sqrt{B(t)}$.
Here the term $J_1$ in (\ref{pin}), reflecting the diffusion mechanism in (\ref{0}) and thus inhibiting
the tendency toward blow-up, can be estimated from above as follows.
\begin{lem}\label{lem52}
  Let $n\ge 1, m>0, K>1$ and $T>0$,
  and suppose that $B\in C^1([0,T)$ is positive and such that (\ref{A_well}) holds as well as
  $K\sqrt{B(t)}<R^n$ for all $t\in [0,T)$.
  Then writing $\xi=\frac{s}{B(t)}$, for the function $J_1$ introduced in (\ref{j1}) we have
  \be{52.1}
	J_1(s,t) \le n B^{1-\frac{1}{n}}(t) \xi^{1-\frac{1}{n}}
	\qquad \mbox{for all $t\in (0,T)$ and } s\in (B(t),K\sqrt{B(t)}).
  \ee
\end{lem}
\proof
  Since $\xi>1$ and hence $\varphi''(\xi)<0$ by (\ref{phi_pp}), we have $|\varphi''(\xi)|=-\varphi''(\xi)$, so that
  we may use the trivial estimate
  \bas
	B^{\frac{4}{n}-2}(t) \varphi'^2(\xi)
	+ n^2 B^{\frac{2}{n}-2}(t) \xi^{2-\frac{2}{n}} \varphi''^2(\xi)
	\ge n^2 B^{\frac{2}{n}-2}(t) \xi^{2-\frac{2}{n}} \varphi''^2(\xi)
  \eas
  to infer that
  \bas
	J_1(s,t)
	&=& n^2 \cdot \frac{\xi^{2-\frac{2}{n}} |\varphi''(\xi)|}
	{\sqrt{B^{\frac{4}{n}-2}(t) \varphi'^2(\xi) + n^2 B^{\frac{2}{n}-2}(t) \xi^{2-\frac{2}{n}} \varphi''^2(\xi)}} \\
	&\le& n^2 \cdot \frac{\xi^{2-\frac{2}{n}} |\varphi''(\xi)|}
		{\sqrt{n^2 B^{\frac{2}{n}-2}(t) \xi^{2-\frac{2}{n}} \varphi''^2(\xi)}} \\
	&=& n B^{1-\frac{1}{n}}(t)\xi^{1-\frac{1}{n}}
  \eas
  holds for any such $t$ and $s$, as claimed.
\qed
Our goal will accordingly consist of controlling the term $J_1$ in (\ref{pin}) from above by a suitably negative
quantity.
As a first step toward this, we shall make sure that
in the root appearing in the denominator of (\ref{j2}), the second summand essentially dominates
the first upon appropriate choices of the parameters.
\begin{lem}\label{lem50}
  Let $n\ge 1, m>0, \lambda\in (0,1), K>1$ with $K\ge \sqrt{\bl R^n}$ and $B_0\in (0,1)$ be such that
  $K\sqrt{B_0}<R^n$ and
  \be{50.01}
	B_0 \le \frac{K^2}{4(\al+\bl)^2}.
%	B_0 \le \frac{K^2}{\bl^2}.
  \ee
  Suppose that for some $T>0$, $B\in C^1([0,T))$ is positive and nonincreasing and such that $B(0)\le B_0$.
  Then writing $\xi=\frac{s}{B(t)}$ for $s\ge 0$ and $t\ge 0$, we have
  \be{50.1}
	\frac{1}{A^2(t) B^{\frac{2}{n}-2}(t) \xi^{\frac{2}{n}-2} \varphi^2(\xi)}
	\le \frac{\omega_n^2}{\lambda^2 m^2} \cdot \Big( 1+\frac{\al R^n}{K^2}\Big) \cdot K^{2-\frac{2}{n}}
	B_0^{3-\frac{3}{n}}
	\qquad \mbox{for all $t\in (0,T)$ and } s\in (B(t),K\sqrt{B(t)}).
  \ee
\end{lem}
\proof
  Since $K\ge \sqrt{\bl R^n}$ and (\ref{50.01}) holds, we know from Lemma \ref{lem44} that $A(t)\ge A_T$ for all
  $t\in (0,T)$ with $A_T$ given by (\ref{AT}).
  Moreover, the fact that $\varphi$ is increasing on $[1,\infty)$ allows us to estimate $\varphi(\xi)\ge 1$
  for all $t\in (0,T)$ and $s>B(t)$, because for any such $t$ and $s$ we have $\xi>1$. Hence,
  \be{50.2}
	\frac{1}{A^2(t) B^{\frac{2}{n}-2}(t) \xi^{\frac{2}{n}-2} \varphi^2(\xi)}
	\le \frac{1}{\lambda^2 A_T^2} \cdot B^{2-\frac{2}{n}}(t) \xi^{2-\frac{2}{n}}
	\qquad \mbox{for all $t\in (0,T)$ and } s\in (B(t),K\sqrt{B(t)}).
  \ee
  As $2-\frac{2}{n}\ge 0$, we may use the restriction $\xi< \frac{K}{\sqrt{B(t)}}$ implied by the inequality
  $s<K\sqrt{B(t)}$ to estimate
  \bas
	\xi^{2-\frac{2}{n}} \le K^{2-\frac{2}{n}} B^{1-\frac{1}{n}}(t)
	\qquad \mbox{for all $t\in (0,T)$ and } s\in (B(t),K\sqrt{B(t)}).
  \eas
  Therefore, (\ref{50.1}) is a consequence of (\ref{50.2}).
\qed
In order to prepare an estimate for the numerator in (\ref{j2}) from below,
let us state and prove the following elementary calculus lemma.
\begin{lem}\label{lem51}
  For $\lambda\in (0,1)$, let $\al$ and $\bl$ be as defined in (\ref{abl}), and let
  \be{51.1}
	\psil(\xi):=\frac{\xi(\xi-\bl)}{\xi-\al-\bl}
	\qquad \mbox{for } \xi\ge 1.
  \ee
  Then if the numbers $K>1$ and $B\in (0,1)$ satisfy
  \be{51.2}
	B \le \frac{K^2}{4(\al+\bl)^2},
  \ee
  we have
  \be{51.3}
	\psil(\xi) \le \max \Big\{ \frac{1}{\lambda} \, , \, \frac{2K}{\sqrt{B}} \Big\}
	\qquad \mbox{for all } \xi \in \Big[1,\frac{K}{\sqrt{B}} \Big].
  \ee
\end{lem}
\proof
  Differentiation in (\ref{51.1}) yields
  \bas
	\psil'(\xi)
	&=& \frac{(2\xi-\bl)(\xi-\al-\bl) - (\xi^2-\bl \xi)}{(\xi-\al-\bl)^2} \\
	&=& \frac{\xi^2 - 2(\al+\bl) \xi + (\al+\bl)\bl}{(\xi-\al-\bl)^2}
	\qquad \mbox{for all } \xi>1,
  \eas
  from which we obtain that
  \be{51.4}
	\psil'(\xi)<0
	\qquad \mbox{if and only if} \qquad
	\xi \in (\xi_-,\xi_+),
  \ee
  where $\xi_+$ and $\xi_-$ are given by
  \bas
	\xi_\pm = \al+\bl \pm \sqrt{(\al+\bl)^2 - (\al+\bl)\bl}.
  \eas
  Here by (\ref{abl}), we recally that $\al+\bl=\frac{\lambda+1}{2}$ in computing
  \bas
	\xi_\pm
	&=&\frac{\lambda+1}{2} \pm \sqrt{\Big(\frac{\lambda+1}{2}\Big)^2 - \frac{\lambda+1}{2} \cdot
		\frac{3\lambda-1}{2\lambda}} \\
	&=& \frac{\lambda+1}{2} \pm \sqrt{\frac{(\lambda+1) \cdot [\lambda(\lambda+1)-(3\lambda-1)]}{4\lambda}} \\
	&=& \frac{\lambda+1}{2} \pm \sqrt{\frac{(\lambda+1) (1-\lambda)^2}{4\lambda}} \\
	&=& \frac{1}{2} \cdot \Big\{ \lambda+1 \pm (1-\lambda) \cdot \sqrt{\frac{\lambda+1}{\lambda}} \Big\}.
  \eas
  Hence,
  \bas
	2(\xi_\pm -1 )
	&=& \lambda-1 \pm (1-\lambda)\cdot \sqrt{\frac{\lambda+1}{\lambda}} \\
	&=& (1-\lambda) \cdot \Big( - 1 \pm \sqrt{\frac{\lambda+1}{\lambda}}\Big),
  \eas
  implying that $\xi_-<1<\xi_+$.
  Therefore, (\ref{51.4}) entails that
  \be{51.5}
	\psil(\xi) \le \max \Big\{ \psil(1) \, , \, \psil\Big(\frac{K}{\sqrt{B}}\Big) \Big\}
	\qquad \mbox{for all } \xi \in \Big[1,\frac{K}{\sqrt{B}} \Big],
  \ee
  where
  \bas
	\psil(1)
	= \frac{1-\bl}{1-\al-\bl}
	= \frac{1-\frac{3\lambda-1}{2\lambda}}{1-\frac{\lambda+1}{2}}
	= \frac{\frac{1-\lambda}{2\lambda}}{\frac{1-\lambda}{2}}
	=\frac{1}{\lambda}.
  \eas
  Since (\ref{51.2}) ensures that $\al+\bl \le \frac{K}{2\sqrt{B}}$ and thus
  \bas
	\psil\Big(\frac{K}{\sqrt{B}}\Big)
	&=& \frac{\frac{K^2}{B}-\bl \cdot \frac{K}{\sqrt{B}}}{\frac{K}{\sqrt{B}}-\al-\bl} \\
	&\le& \frac{\frac{K^2}{B}}{\frac{K}{\sqrt{B}} - \al-\bl} \\
	&\le& \frac{\frac{K^2}{B}}{\frac{K}{2\sqrt{B}}} \\
	&=& \frac{2K}{\sqrt{B}},
  \eas
  the inequality (\ref{51.5}) thus yields (\ref{51.3}).
\qed
On the basis of the above lemma, we can indeed achieve that in the numerator in (\ref{j2}) the positive
summand prevails.
\begin{lem}\label{lem53}
  Let $n\ge 1, m>0, \lambda\in (0,1), K>1, \delta\in (0,1)$ and $B_0\in (0,1)$ such that $K\sqrt{B_0}<R^n$ and
  \be{53.1}
	B_0 \le \frac{K^2}{4(\al+\bl)^2}
  \ee
  as well as
  \be{53.2}
	\frac{\mu}{n A_T} \cdot \max \Big\{ \frac{B_0}{\lambda} \, , \, 2K\sqrt{B_0} \Big\} \le \delta
  \ee
  with $\mu$ and $A_T$ as in (\ref{mu}) and (\ref{AT}), respectively.
  Furthermore, let $T>0$ and $B\in C^1([0,T))$ be positive and such that
  \be{53.22}
	B(t) \le B_0
	\qquad \mbox{for all } t\in (0,T).
  \ee
  Then writing $\xi=\frac{s}{B(t)}$, we have
  \be{53.3}
	A(t)\varphi(\xi) - \frac{\mu}{n} B(t)\xi \ge (1-\delta) A(t)\varphi(\xi)
	\qquad \mbox{for all $t\in (0,T)$ and } s\in (B(t),K\sqrt{B(t)}).
  \ee
\end{lem}
\proof
  With $\psil$ taken from Lemma \ref{lem51}, we first observe that
  \bea{53.4}
	A(t)\varphi(\xi)-\frac{\mu}{n} B(t)\xi
	&=& A(t)\varphi(\xi) \cdot \Big\{ 1 - \frac{\mu B(t) \xi}{nA(t)} \cdot \frac{\xi-\bl}{\xi-\al-\bl}\Big\}
		\nn\\
	&=& A(t)\varphi(\xi) \cdot \Big\{ 1 - \frac{\mu B(t)}{nA(t)} \cdot \psil(\xi) \Big\}
%	\qquad \mbox{for all $t\in (0,T)$ and } s\in (B(t),K\sqrt{B(t)}).
  \eea
  for all $t\in (0,T)$ and $s\in (B(t),K\sqrt{B(t)})$.
  Here thanks to (\ref{53.1}) and (\ref{53.22}) we may apply Lemma \ref{lem51}, which combined with
  (\ref{AT}) shows that
  \bas
	\frac{\mu B(t)}{nA(t)} \cdot \psil(\xi)
	&\le& \frac{\mu B(t)}{nA_T} \cdot \max \Big\{ \frac{1}{\lambda} \, , \, \frac{2K}{\sqrt{B(t)}} \Big\} \\
	&\le& \frac{\mu}{nA_t} \cdot \max \Big\{ \frac{B_0}{\lambda} \, , \, 2K\sqrt{B_0} \Big\}
	\qquad \mbox{for all $t\in (0,T)$ and } s\in (B(t),K\sqrt{B(t)}).
  \eas
  In light of (\ref{53.2}), the conclusion (\ref{53.3}) is therefore a consequence of (\ref{53.4}).
\qed
With the above preparations at hand, we can proceed to show that under the assumptions of Theorem \ref{theo14},
if $B$ is a suitably small nonincreasing function
satisfying an appropriate differential inequality, then $\win$ indeed becomes a subsolution of (\ref{0w})
in the intermediate region where $B(t)<s<K\sqrt{B(t)}$.\\
We shall first demonstrate this in the spatially one-dimensional case, in which the role of the number $m_c$ in
(\ref{mc}) will become clear through the following lemma.
\begin{lem}\label{lem54}
  Let $n=1, \chi>1$ and $m>m_c=\frac{1}{\sqrt{\chi^2-1}}$.
  Then there exist $\lambda\in (0,1)$, $\kappa>0$ and $B_0\in (0,1)$ such that $K\sqrt{B_0}<R$, and such that
  whenever $T>0$ and $B\in C^1([0,T))$ is a positive and nonincreasing function fulfilling (\ref{A_well}) as well as
  \be{54.1}
	\left\{ \begin{array}{l}
	B'(t) \ge - \kappa \sqrt{B(t)}, \qquad t\in (0,T), \\[1mm]
	B(0) \le B_0,
	\end{array} \right.
  \ee
  then for $\win$ as in (\ref{win}) we have
  \be{54.2}
	(\parab \win)(s,t)\le 0
	\qquad \mbox{for all $t\in (0,T)$ and } s\in (B(t),K\sqrt{B(t)}).
  \ee
\end{lem}
\proof
  As $m>m_c$, we have $\frac{m\chi}{\sqrt{1+m^2}}>1$, whence it is possible to fix $\lambda\in (0,1)$ sufficiently
  close to $1$ such that
  \bas
	\frac{m\chi}{\sqrt{\frac{1}{\lambda^2}+m^2}} >1.
  \eas
  This in turn allows us to choose some $\delta\in (0,1)$ such that
  \be{54.3}
	c_1:=\frac{(1-\delta)m\chi}{\sqrt{\frac{1+\delta}{\lambda^2}+m^2}} -1
  \ee
  is positive.
  We thereafter pick $K>1$ such that with $\al$ and $bl$ as in (\ref{abl}) we have
  \be{54.4}
	K\ge \sqrt{\bl R}
  \ee
  and
  \be{54.5}
	\frac{\al R}{K^2} \le \delta.
  \ee
  Finally, we take $B_0\in (0,1)$ conveniently small fulfilling $K\sqrt{B_0}<R$ and
  \be{54.6}
	B_0 \le \frac{K^2}{4(\al+\bl)^2}
  \ee
  as well as
  \be{54.7}
	\frac{\mu}{A_T} \cdot \max \Big\{ \frac{B_0}{\lambda} \, , \, 2K\sqrt{B_0} \Big\} \le \delta
  \ee
  with $A_T$ as in (\ref{AT}), and let
  \be{54.8}
	\kappa:=\frac{c_1}{K}.
  \ee
  Then given any $T>0$ and a positive nonincreasing $B\in C^1([0,T))$ satisfying (\ref{54.1}), from Lemma
  \ref{lem44} in conjunction with (\ref{54.4}) we know that $A' \le 0$ on $(0,T)$, so that
  (\ref{pin}) yields
  \be{54.80}
	\frac{B(t)}{A(t)\varphi'(\xi)} \cdot (\parab \win)(s,t)
	\le -\xi B'(t) + J_1(s,t)+J_2(s,t)
	\qquad \mbox{for all $t\in (0,T)$ and } s\in (B(t), K\sqrt{B(t)})
  \ee
  with $\xi=\frac{s}{B(t)}$ and $J_1$ and $J_2$ as given by (\ref{j1}) and (\ref{j2}).\\
  Here, Lemma \ref{lem52} says that
  \be{54.9}
	J_1(s,t) \le 1
	\qquad \mbox{for all $t\in (0,T)$ and } s\in (B(t), K\sqrt{B(t)}),
  \ee
  and in order to compensate this positive contribution in (\ref{54.80}) appropriately, we first invoke
  Lemma \ref{lem53}, which ensures that thanks to (\ref{54.6}) and (\ref{54.7}) we have
  \be{54.10}
	A(t)-\mu B(t)\xi \ge (1-\delta)A(t)\varphi(\xi)
	\qquad \mbox{for all $t\in (0,T)$ and } s\in (B(t), K\sqrt{B(t)}).
  \ee
  In particular, this implies that the expression on the left-hand side herein is nonnegative, so that we can
  estimate
  \be{54.11}
	\Big(A(t)\varphi(\xi)-\mu B(t)\xi\Big)^2 \le A^2(t)\varphi^2(\xi)
	\qquad \mbox{for all $t\in (0,T)$ and } s\in (B(t), K\sqrt{B(t)}).
  \ee
  Since (\ref{54.4}) and (\ref{54.6}) allow for an application of Lemma \ref{lem50}, we moreover know that
  \bas
	\frac{1}{A^2(t)\varphi^2(\xi)}
	&\le& \frac{1}{\lambda^2 m^2} \cdot \Big(1+\frac{\al R}{K^2}\Big) \\
	&\le& \frac{1+\delta}{\lambda^2 m^2}
	\qquad \mbox{for all $t\in (0,T)$ and } s\in (B(t), K\sqrt{B(t)})
  \eas
  because of (\ref{54.5}).
  Combining this with (\ref{54.11}) shows that in the denominator in the definition (\ref{j2}) of $J_2$ we have
  \bas
	\sqrt{1+\Big(A(t)\varphi(\xi)-\mu B(t)\xi\Big)^2}
	&\le& \sqrt{\frac{1+\delta}{\lambda^2 m^2} \cdot A^2(t)\varphi^2(\xi) + A^2(t)\varphi^2(\xi)} \\
	&=& \sqrt{\frac{1+\delta}{\lambda^2 m^2}+1} \cdot A(t)\varphi(\xi)
	\qquad \mbox{for all $t\in (0,T)$ and } s\in (B(t), K\sqrt{B(t)}),
  \eas
  so that by means of (\ref{54.10}) we can estimate
  \bas
	-J_2(s,t)
	&=& \chi \cdot \frac{A(t)\varphi(\xi)-\mu B(t)\xi}
	{1+\Big(A(t)\varphi(\xi)-\mu B(t)\xi\Big)^2} \\
	&\ge& \chi \cdot \frac{(1-\delta) A(t)\varphi(\xi)}{\sqrt{\frac{1+\delta}{\lambda^2 m^2} +1} \cdot
		A(t)\varphi(\xi)} \\
	&=& \frac{(1-\delta) m\chi}{\sqrt{\frac{1+\delta}{\lambda^2} + m^2}}
	\qquad \mbox{for all $t\in (0,T)$ and } s\in (B(t), K\sqrt{B(t)}).
  \eas
  Together with (\ref{54.8}) and (\ref{54.9}), in view of the definition (\ref{54.3}) of $c_1$ this implies that
  \bas
	\frac{B(t)}{A(t)\varphi'(\xi)} \cdot (\parab \win)(s,t)
	&\le& -\xi B'(t) + 1 - \frac{(1-\delta)m\chi}{\sqrt{\frac{1+\delta}{\lambda^2}+m^2}} \\[2mm]
	&=& -\xi B'(t) - c_1
	\qquad \mbox{for all $t\in (0,T)$ and } s\in (B(t), K\sqrt{B(t)}).
  \eas
  Once more using that in the considered region we have $\xi\le \frac{K}{\sqrt{B(t)}}$, due to our choice of
  $\kappa$ we infer that
  \bas
	-\xi B'(t) + c_1
	&=& \xi \cdot \Big\{ -B'(t)-\frac{c_1}{\xi}\Big\} \\
	&\le& \xi \cdot \Big\{ -B'(t) - \frac{c_1 \sqrt{B(t)}}{K}\Big\} \\
	&=& \xi \cdot \Big\{ -B'(t) - \kappa\sqrt{B(t)} \Big\} \\[2mm]
	&\le& 0
	\qquad \mbox{for all $t\in (0,T)$ and } s\in (B(t), K\sqrt{B(t)})
  \eas
  because of (\ref{54.1}), whereby the proof is completed.
\qed
In the case $n\ge 2$, we follow the same basic strategy as above, but numerous adaptations are necessary due to the fact
that in this case the more involved, and more degenerate, structure of $J_1$ and $J_2$ in (\ref{pin})
allow for choosing actually any positive value of the mass $m$ whenever $\chi>1$.
\begin{lem}\label{lem55}
  Let $n\ge 2, \chi>1$ and $m>0$, and let $\lambda\in (0,1)$ be arbitrary.
  Then there exist $K>1$, $\kappa>0$ and $B_0\in (0,1)$ such that $K\sqrt{B_0}<R^n$, and such that
  if $T>0$ and $B\in C^1([0,T))$ is positive and nonincreasing such that
  \be{55.1}
	\left\{ \begin{array}{l}
	B'(t) \ge - \kappa B^{1-\frac{1}{2n}}(t), \qquad t\in (0,T), \\[1mm]
	B(0) \le B_0,
	\end{array} \right.
  \ee
  then the function $\win$ defined in (\ref{wout}) satisfies
  \be{55.2}
	(\parab \win)(s,t)\le 0
	\qquad \mbox{for all $t\in (0,T)$ and } s\in (B(t),K\sqrt{B(t)}).
  \ee
\end{lem}
\proof
  We let $\al$ and $\bl$ as in (\ref{abl}), take any $K>1$ fulfilling
  \be{55.22}
	K>\sqrt{\bl R^n}
  \ee
  and use that $\chi>1$ to pick $\delta\in (0,1)$ suitably small such that
  \be{55.3}
	c_1:=n\cdot \bigg\{ \frac{(1-\delta)\chi}{\sqrt{1+\delta}} - 1 \bigg\} \ >0.
  \ee
  It is the possible to fix $B_0\in (0,1)$ such that $K\sqrt{B_0}<R^n$ and
  \be{55.4}
	B_0 \le \frac{K^2}{4(\al+\bl)^2},
  \ee
  such that with $A_T$ as in (\ref{AT}) we have
  \be{55.5}
	\frac{\mu}{nA_t} \cdot \max \Big\{ \frac{B_0}{\lambda} \, , \, 2K\sqrt{B_0} \Big\} \le \delta,
  \ee
  and such that
%  \be{55.6}
%	B_0 le \frac{K}{\bl}
%  \ee
%  as well as
  \be{55.7}
	\frac{\omega_n}{\lambda^2 m^2} \cdot \Big( 1+\frac{\al R^n}{K^2}\Big) \cdot K^{2-\frac{2}{n}}
	B_0^{3-\frac{3}{n}} \le \delta,
  \ee
  where we not that in achieving the latter we make use of our assumption that $n\ge 2$.
  We finally let
  \be{55.8}
	\kappa:= c_1 K^{-\frac{1}{n}},
  \ee
  and suppose that $T>0$ and that $B\in C^1([0,T))$ is positive and nonincreasingand such that (\ref{55.1}) holds.\\
  Then (\ref{55.22}) and (\ref{55.4}) warrant that Lemma \ref{lem44} applies so as to yield that
  $A'\le 0$ on $(0,T)$, and that hence by (\ref{pin}),
  \be{55.9}
	\frac{B(t)}{A(t)\varphi'(\xi)} \cdot (\parab \win)(s,t)
	\le -\xi B'(t) + J_1(s,t)+J_2(s,t)
	\qquad \mbox{for all $t\in (0,T)$ and } s\in (B(t),K\sqrt{B(t)}),
  \ee
  where again $\xi=\frac{s}{B(t)}$, and where $J_1$ and $J_2$ are as defined in (\ref{j1}) and (\ref{j2}),
  respectively.
  Now thanks to (\ref{55.4}) and (\ref{55.5}), Lemma \ref{lem53} shows that
  \be{55.10}
	A(t)\varphi(\xi) - \frac{\mu}{n} B(t)\xi
	\ge (1-\delta) A(t)\varphi(\xi)
	\qquad \mbox{for all $t\in (0,T)$ and } s\in (B(t),K\sqrt{B(t)}),
  \ee
  whereas (\ref{55.4}) allows for invoking Lemma \ref{lem50} to infer from (\ref{55.7}) that
  \bea{55.11}
	\frac{1}{A^2(t) B^{\frac{2}{n}-2}(t) \xi^{\frac{2}{n}-2} \varphi^2(\xi)}
	&\le& \frac{\omega_n^2}{\lambda^2 m^2} \cdot \Big(1+\frac{\al R^n}{K^2}\Big) \cdot K^{2-\frac{2}{n}}
		B_0^{3-\frac{3}{n}}  \nn\\[2mm]
	&\le& \delta
	\qquad \mbox{for all $t\in (0,T)$ and } s\in (B(t),K\sqrt{B(t)}).
  \eea
  By means of (\ref{55.10}), (\ref{55.11}) and the fact that $\delta<1$, we can thus estimate $J_2$ according to
  \bas
	-J_2(s,t)
	&=& n\chi \cdot \frac{A(t)\varphi(\xi)-\frac{\mu}{n}B(t)\xi}
	{\sqrt{1+B^{\frac{2}{n}-2}(t) \xi^{\frac{2}{n}-2} \Big(A(t)\varphi(\xi)-\frac{\mu}{n}B(t)\xi \Big)^2}} \nn\\
	&\ge& n\chi \cdot \frac{(1-\delta) A(t)\varphi(\xi)}
	{\sqrt{1+B^{\frac{2}{n}-2}(t) \xi^{\frac{2}{n}-2} \Big(A(t)\varphi(\xi)-\frac{\mu}{n}B(t)\xi \Big)^2}} \nn\\
	&\ge& n\chi \cdot \frac{(1-\delta)A(t)\varphi(\xi)}
	{\sqrt{1+B^{\frac{2}{n}-2}(t) \xi^{\frac{2}{n}-2} A^2(t) \varphi^2(\xi)}} \\
	&\ge& n\chi \cdot \frac{(1-\delta)A(t)\varphi(\xi)}
	{\sqrt{(\delta+1) \cdot B^{\frac{2}{n}-2}(t) \xi^{\frac{2}{n}-2} A^2(t) \varphi^2(\xi)}} \\[1mm]
	&=& \frac{(1-\delta) n\chi}{\sqrt{1+\delta}} \cdot B^{1-\frac{1}{n}}(t) \xi^{1-\frac{1}{n}}
	\qquad \mbox{for all $t\in (0,T)$ and } s\in (B(t),K\sqrt{B(t)}).
  \eas
  Since on the other hand
  \bas
	J_1(s,t) \le n B^{1-\frac{1}{n}}(t) \xi^{1-\frac{1}{n}}
	\qquad \mbox{for all $t\in (0,T)$ and } s\in (B(t),K\sqrt{B(t)})
  \eas
  due to Lemma \ref{lem53}, we therefore conclude from (\ref{55.9}) that
  \bea{55.12}
	\frac{B(t)}{A(t)\varphi'(\xi)} \cdot (\parab \win)(s,t)
	&\le& -\xi B'(t) + n B^{1-\frac{1}{n}}(t) \xi^{1-\frac{1}{n}}
	- \frac{(1-\delta)n\chi}{\sqrt{1+\delta}} \cdot B^{1-\frac{1}{n}}(t) \xi^{1-\frac{1}{n}} \nn\\
	&=& - \xi B'(t) - c_1 B^{1-\frac{1}{n}}(t) \xi^{1-\frac{1}{n}}
%	\qquad \mbox{for all $t\in (0,T)$ and } s\in (B(t),K\sqrt{B(t)}).
  \eea
  for all $t\in (0,T)$ and $s\in (B(t), K\sqrt{B(t)})$.
  We finally observe that $\xi<\frac{K}{\sqrt{B(t)}}$ whenever $s<K\sqrt{B(t)}$, and that hence by (\ref{55.8}),
  \bas
	-\xi B'(t) - c_1 B^{1-\frac{1}{n}}(t) \xi^{1-\frac{1}{n}}
	&=& \xi \cdot \Big\{ - B'(t) - c_1 B^{1-\frac{1}{n}}(t) \xi^{-\frac{1}{n}} \Big\} \\
	&\le& \xi \cdot \Big\{ -B'(t) - c_1 B^{1-\frac{1}{n}}(t) \cdot K^{-\frac{1}{n}} B^\frac{1}{2n}(t)\Big\} \\
	&=& \xi \cdot \Big\{ -B'(t) - \kappa B^{1-\frac{1}{2n}}(t)\Big\}
%	\qquad \mbox{for all $t\in (0,T)$ and } s\in (B(t),K\sqrt{B(t)})
  \eas
  for all $t\in (0,T)$ and $s\in (B(t), K\sqrt{B(t)})$,
  so that (\ref{55.1}) and (\ref{55.12}) guarantee that indeed the claimed inequality (\ref{55.2}) holds.
\qed
\mysection{Blow-up. Proof of Theorem \ref{theo14}}\label{sect4}
Now our final result on blow-up of solutions to the original problem
can be derived by a combination of Lemma \ref{lem60} with Lemma \ref{lem61} as well Lemma \ref{lem54} and
Lemma \ref{lem55} in the cases $n=1$ and $n\ge 2$, respectively,
along with a straightforward comparison argument.\abs
\proofc of Theorem \ref{theo14}.\quad
  Thanks to our assumptions (\ref{chi}) and (\ref{m}), in view of Lemma \ref{lem60}, Lemma \ref{lem61},
  Lemma \ref{lem54} and Lemma \ref{lem55} we can fix $\lambda\in (0,1), K>0, \kappa>0$ and $B_0\in (0,1)$ such that
  $K\sqrt{B_0}<R^n$ and such that if we let $B$ denote the solution of
  \be{B}
	\left\{ \begin{array}{l}
	B'(t)=-\kappa B^{1-\frac{1}{2n}}(t), \qquad t\in (0,T), \\[1mm]
	B(0)=B_0,
	\end{array} \right.
  \ee
  extended up to its extinction time $T\in (0,\infty)$, that is, if we define
  \be{B1}
	B(t):=\Big\{ B_0^\frac{1}{2n} - \frac{\kappa}{2n}t \Big\}^{2n}, \qquad t\in [0,T),
  \ee
  with
  \be{T}
	T:=\frac{2n}{\kappa} \cdot B_0^\frac{1}{2n},
  \ee
  then the functions $\wout$ and $\win$ given by (\ref{wout}) and (\ref{win}) are well-defined and satisfy
  \be{14.4}
	(\parab \wout)(s,t) \le 0
	\qquad \mbox{for all $t\in (0,T)$ and } s\in (K\sqrt{B(t)},R^n)
  \ee
  as well as
  \be{14.5}
	(\parab \win)(s,t) \le 0
	\qquad \mbox{for all $t\in (0,T)$ and } s\in (0,B(t)) \cup (B(t),K\sqrt{B(t)}).
  \ee
  Here in applying Lemma \ref{lem61} we note that $-\kappa B^{1-\frac{1}{2n}}(t) \ge -\kappa B^{1-\frac{1}{n}}(t)$
  for all $t\in (0,T)$ due to the fact that $B(t)\le B_0<1$ for all $t\in (0,T)$.
  According to (\ref{14.4}) and (\ref{14.5}), Lemma \ref{lem21} asserts that
  \bas
	\uw(s,t):=\left\{ \begin{array}{ll}
	\win(s,t) \qquad & \mbox{if $t\in [0,T)$ and } s\in [0,K\sqrt{B(t)}], \\[1mm]
	\wout(s,t) \qquad & \mbox{if $t\in [0,T)$ and } s\in (K\sqrt{B(t)}, R^n],
	\end{array} \right.
  \eas
  defines a function $\uw\in C^1([0,R^n]\times [0,T))$ which satisfies
  $\uw(\cdot,t)\in C^2([0,R^n] \setminus \{B(t), K\sqrt{B(t)}\})$ for all $t\in [0,T)$ as well as
  \bas
	(\parab \uw)(s,t) \le 0
	\qquad \mbox{for all $t\in [0,T)$ and } s\in (0,R^n) \setminus \{B(t),K\sqrt{B(t)}\}.
  \eas
  Therefore, if $u_0$ satisfies (\ref{init}) and is such that
  \bas
	\int_{B_r(0)} u_0(x)dx \ge M(r) := \omega_n \uw(r^n,0)
	\qquad \mbox{for all } r\in [0,R],
  \eas
  then the solution $w$ of (\ref{0w}) defined through (\ref{trans}) satisfies
  \be{14.6}
	w(s,0)\ge \uw(s,0)
	\qquad \mbox{ for all $s\in (0,R^n)$},
  \ee
  and furthermore it is clear that
  \be{14.7}
	w(0,t)=\uw(0,t)=0
	\quad \mbox{and} \quad
	w(R^n,t)=\uw(R^n,t)=\frac{m}{\omega_n}
	\qquad \mbox{ for all } t\in (0,\tilde T),
  \ee
  where $\tilde T:=\min \{\tm, T\}$.
  In order to assert applicability of
  the comparison principle from Lemma \ref{lem70} below, we abbreviate $\alpha:=2-\frac{2}{n}\ge 0$ and let
  \bas
	\phi(s,t,y_0,y_1,y_2)
	:=n^2 \cdot \frac{s^\alpha y_1 y_2}{\sqrt{y_1^2 + n^2 s^\alpha y_2^2}}
	+ n\chi \cdot \frac{(y_0-\frac{\mu}{n}s) y_1}{\sqrt{1+s^{-\alpha} (y_0-\frac{\mu}{n}s)^2}}
%	\qquad \mbox{for } (s,t,y_0,y_1,y_2) \in G:=(0,R^n)\times (0,\infty)\times \R \times (0,\infty) \times \R,
  \eas
  for $(s,t,y_0,y_1,y_2) \in G:=(0,R^n)\times (0,\infty)\times \R \times (0,\infty) \times \R$,
  so that $\phi\in C^1(G)$ with
  \bea{14.8}
	\frac{\partial\phi}{\partial y_2}(s,t,y_0,y_1,y_2)
	&=& n^2 \cdot \frac{s^\alpha y_1^3}{\sqrt{y_1^2+n^2 s^\alpha y_2^2}} \nn\\[1mm]
	&\ge& 0
	\qquad \mbox{for all } (s,t,y_0,y_1,y_2)\in G
  \eea
  and
  \bas
	\frac{\partial\phi}{\partial y_1}(s,t,y_0,y_1,y_2)
	= n^4 \cdot \frac{s^{2\alpha} y_2^3}{\sqrt{y_1^2 + n^2 s^\alpha y_2^2}}
	+ n\chi \cdot \frac{y_0-\frac{\mu}{n}s}{\sqrt{1+s^{-\alpha} (y_0-\frac{\mu}{n}s)^2}}
	\qquad \mbox{for all } (s,t,y_0,y_1,y_2)\in G
  \eas
  as well as
  \bas
	\frac{\partial\phi}{\partial y_0}(s,t,y_0,y_1,y_2)
	= n\chi \cdot \frac{y_1}{\sqrt{1+s^{-\alpha} (y_0-\frac{\mu}{n}s)^2}}
	\qquad \mbox{for all } (s,t,y_0,y_1,y_2)\in G.
  \eas
  Therefore, we can estimate
  \bea{14.9}
	\Big|
	\frac{\partial\phi}{\partial y_1}(s,t,y_0,y_1,y_2)
	\Big|
	&\le& ns^\frac{\alpha}{2} \cdot
	\frac{\sqrt{n^2 s^\alpha y_2^2}^3}{\sqrt{y_1^2+ n^2 s^\alpha y_2^2}^3}
	+ n\chi s^\frac{\alpha}{2} \cdot
	\frac{\sqrt{s^{-\alpha}(y_0-\frac{\mu}{n}s)^2}}{\sqrt{1+s^{-\alpha} (y_0-\frac{\mu}{n}s)^2}} \nn\\
	&\le& n R^\frac{2\alpha}{n} + n\chi R^\frac{2\alpha}{n}
	\qquad \mbox{for all } (s,t,y_0,y_1,y_2)\in G
  \eea
  and
  \be{14.10}
	\Big|
	\frac{\partial\phi}{\partial y_0}(s,t,y_0,y_1,y_2)
	\Big|
	\le n\chi |y_1|
	\qquad \mbox{for all } (s,t,y_0,y_1,y_2)\in G.
  \ee
  Since the inequalities (\ref{14.8}), (\ref{14.9}) and (\ref{14.10}) warrant the validity of the hypotheses
  (\ref{70.1}), (\ref{70.3}) and (\ref{70.2}) of Lemma \ref{lem70}, as a consequence of the latter we obtain that
  \bas
	w(s,t) \ge \uw(s,t)
	\qquad \mbox{for all $s\in [0,R^n]$ and } t\in [0,\tilde T).
  \eas
  As $w(0,t)=\uw(0,t)=0$ for all $t\in (0,\tilde T)$, by the mean value theorem this implies that
  for each $t\in (0,\tilde T)$ we can find some $\theta(t)\in (0,R^n)$ with the property that
  \bas
	w_s(\theta(t),t)
	= \frac{w(B(t),t)}{B(t)}
	\ge \frac{\uw(B(t),t)}{B(t)}
	= \frac{A(t)\varphi(1)}{B(t)}
	= \lambda \cdot \frac{A(t)}{B(t)}
	\qquad \mbox{for all } t\in (0,\tilde T).
  \eas
  Recalling that $u(r,t)=w_s(r^\frac{1}{n},t)$ for all $r\in (0,R)$ and $t\in (0,\tm)$, we thereby infer that
  \bas
	\sup_{r\in (0,R)} u(r,t) \ge w_s(\theta(t),t)=\lambda\cdot \frac{A(t)}{B(t)}
	\qquad \mbox{for all } t\in (0,\tilde T).
  \eas
  In view of the fact that $B(t)\to 0$ as $t\nearrow T$, and that hence $A(t)\to \frac{m}{\omega_n}$ as $t\nearrow T$
  according to (\ref{A}), this entails that we necessarily must have $\tm \le T<\infty$,
  so that (\ref{14.2}) becomes a consequence of the extensibility criterion (\ref{43.1}).
\qed
\mysection{Appendix: A comparison lemma}
An ingredient essential to our argument is the following variant of the parabolic comparison principle.
Since we could not find an appropriate reference precisely covering the present situation, especially involving
the present particular type of degenerate diffusion and nonsmooth comparison functions, we include a proof
for completeness.
\begin{lem}\label{lem70}
  Let $L>0, T>0, G:=(0,L)\times (0,T) \times \R \times (0,\infty)\times \R$ and $\phi \in C^1(G)$ be such that
  \be{70.1}
	\frac{\partial\phi}{\partial y_2} (s,t,y_0,y_1,y_2)\ge 0
	\qquad \mbox{for all } (s,t,y_0,y_1,y_2) \in G,
  \ee
  that for all $T_0\in (0,T)$ and $\Lambda>0$ there exists $C(T_0,\Lambda)>0$ fulfilling
  \be{70.2}
	\Big|\frac{\partial\phi}{\partial y_0}(s,t,y_0,y_1,y_2)\Big| \le C(T_0,\Lambda)
	\qquad \mbox{for all } (s,t,y_0,y_1,y_2) \in G
	\quad \mbox{with $t\in (0,T_0)$ and } y_1\in (0,\Lambda),
  \ee
  and such that for any $t_0\in (0,T)$ we have
  \be{70.3}
	\frac{\partial\phi}{\partial y_1} (\cdot,t,\cdot,\cdot,\cdot)
	\in L^\infty_{loc}((0,L)\times \R \times (0,\infty)\times \R).
  \ee
  Suppose that $\uw$ and $ow$ are two functions which belong to $C^1([0,L]\times [0,T))$ and satisfy
  \be{70.4}
	\uw_s(s,t)>0
	\quad \mbox{and} \quad
	\ow(s,t)>0
	\qquad \mbox{for all $s\in (0,L)$ and } t\in (0,T)
  \ee
  as well as
  \be{70.5}
	\uw(\cdot,t) \in W^{2,\infty}_{loc}((0,L))
	\quad \mbox{and} \quad
	\ow(\cdot,t) \in W^{2,\infty}_{loc}((0,L))
	\qquad \mbox{for all } t\in (0,T).
  \ee
  If moreover
  \be{70.6}
	\uw_t \le \phi(s,t,\uw,\uw_s,\uw_ss)
	\quad \mbox{and}
	\ow_t \ge \phi(s,t,\uw,\uw_s,\uw_ss)
	\qquad \mbox{for all $t\in (0,T)$ and a.e.~$s\in (0,L)$}
  \ee
  and
  \be{70.7}
	\uw(s,0)\le \ow(s,0)
	\qquad \mbox{for all } s\in (0,L)
  \ee
  as well as
  \be{70.8}
	\uw(0,t)\le \ow(0,t)
	\quad \mbox{and} \quad
	\uw(L,t)\le \ow(L,t)
	\qquad \mbox{for all } t\in (0,T),
  \ee
  then
  \be{70.9}
	\uw(s,t)\le \ow(s,t)
	\qquad \mbox{for all $s\in [0,L]$ and } t\in [0,T).
  \ee
\end{lem}
\proof
  We fix an arbitrary $T_0\in (0,T)$ and then obtain from (\ref{70.4}) and the assumed regularity properties of
  $\uw$ and $\ow$ that there exists $\Lambda=\Lambda(T_0)>0$ such that
  \be{70.10}
	0<\uw_s(s,t)<\Lambda
	\quad \mbox{and} \quad
	0<\ow_s(s,t)<\Lambda
	\qquad \mbox{for all $s\in (0,L)$ and } t\in (0,T_0).
  \ee
  For $\eps>0$, we then let $c_1:=C(T_0,\Lambda)$ with $C(T_0,\Lambda)>0$ as in (\ref{70.2}), define
  \be{70.11}
	z(s,t):=\uw(s,t)-\ow(s,t)-\eps \, e^{2c_1 t}
	\qquad \mbox{for $s\in [0,L]$ and } t\in [0,T_0],
  \ee
  and claim that
  \be{70.12}
	z(s,t)<0
	\qquad \mbox{for all $s\in [0,L]$ and } t\in [0,T_0).
  \ee
  To verify this, supposing for contradiction that (\ref{70.12}) be false, from (\ref{70.7}) and (\ref{70.8}) we would infer
  the existence of $s_0\in (0,L)$ and $t_0\in (0,T_0)$ such that
  \be{70.13}
	\max_{(s,t)\in [0,L]\times [0,t_0]} z(s,t)=z(s_0,t_0)=0,
  \ee
  in particular implying that
  \be{70.14}
	z_t(s_0,t_0)\ge 0
  \ee
  and
  \be{70.15}
	z_s(s_0,t_0)=0.
  \ee
  Moreover, using (\ref{70.15}) we obtain that $z(\cdot,t_0)\in W^{2,\infty}_{loc}((0,L))$, so that we can find a null set
  $N\subset (0,L)$ such that $z_{ss}(s,t_0)$ exists for all $s\in (0,L)\setminus N$
  and
%, that the inequalities in (\ref{70.6})
%  are valid at $(s,t_0)$ for all $s\in (0,L)\setminus N$ and
  \be{70.16}
	z_s(s,t_0)=\int_{s_0}^s z_{ss}(\sigma,t_0)d\sigma
	\qquad \mbox{for all } s\in [0,L]
  \ee
  according to (\ref{70.15}), where for later use we note that enlarging $N$ if necessary we can furthermore
  achieve that both inequalities in (\ref{70.6})
  are valid at $(s,t_0)$ for all $s\in (0,L)\setminus N$.
  As $z(\cdot,t_0)$ attains its maximum at $s_0$ by (\ref{70.13}), the identity (\ref{70.16}) necessarily requires that
  there exists $(s_j)_{j\in\N}\subset (s_0,L)\setminus N$ such that $s_j\searrow s_0$ as $j\to\infty$ and
  \be{70.17}
	z_{ss}(s_j,t_0) \le 0
	\qquad \mbox{for all } j\in\N,
  \ee
  for otherwise (\ref{70.16}) would imply that $z_s(s,t_0)>0$ for all $s\in (s_0,s_\star)$ with some $s_\star\in (s_0,L)$,
  which would clearly contradict (\ref{70.13}).\\
  Now differentiating (\ref{70.11}), in view of (\ref{70.6}) and our choice of $N$ we see that
  \bas
	z_t
	&=& \uw_t - \ow_t -2c_1 \eps e^{2c_1 t_0} \nn\\
	&\le& \phi(s,t_0,\uw,\uw_s,\uw_{ss}) - \phi(s,t_0,\ow,\ow_s,\ow_{ss}) - 2c_1 \eps e^{2c_1 t_0}
	\qquad \mbox{for all } s\in N,
  \eas
  so that from (\ref{70.17}) we infer that
  \bea{70.18}
	z_t(s_j,t_0)
	&\le& \phi \Big(s_j,t_0,\uw(s_j,t_0), \uw_s(s_j,t_0), \ow_{ss}(s_j,t_0)\Big)
	- \phi \Big(s_j,t_0,\ow(s_j,t_0), \ow_s(s_j,t_0), \ow_{ss}(s_j,t_0)\Big) \nn\\[1mm]
	& & - 2c_1 \eps e^{2c_1 \eps t_0}
	\qquad \mbox{for all } j\in\N.
  \eea
  Here by the mean value theorem we have
  \bea{70.19}
	& & \hspace*{-30mm}
	\phi \Big(s_j,t_0,\uw(s_j,t_0), \uw_s(s_j,t_0), \ow_{ss}(s_j,t_0)\Big)
	- \phi \Big(s_j,t_0,\ow(s_j,t_0), \ow_s(s_j,t_0), \ow_{ss}(s_j,t_0)\Big) \nn\\
	&=& \xi_j \cdot \Big(\uw(s_j,t_0)-\ow(s_j,t_0\Big)
	+ \lambda_j \cdot \Big(\uw_s(s_j,t_0)-\ow_s(s_j,t_0\Big)
  \eea
  with
  \bas
	\xi_j:=
	\int_0^1 \frac{\partial\phi}{\partial y_0}
	\Big(s_j,t_0,
	\ow(s_j,t_0)+\sigma (\uw(s_j,t_0)-\ow(s_j,t_0)),
	\ow_s(s_j,t_0)+\sigma (\uw_s(s_j,t_0)-\ow_s(s_j,t_0)),
	\ow_{ss}(s_j,t_0)\Big) d\sigma
  \eas
  and
  \bas
	\lambda_j:=
	\int_0^1 \frac{\partial\phi}{\partial y_1}
	\Big(s_j,t_0,
	\ow(s_j,t_0)+\sigma (\uw(s_j,t_0)-\ow(s_j,t_0)),
	\ow_s(s_j,t_0)+\sigma (\uw_s(s_j,t_0)-\ow_s(s_j,t_0)),
	\ow_{ss}(s_j,t_0)\Big) d\sigma
  \eas
  for $j\in\N$.
  Since $s_j\to s_0$ as $j\to\infty$, by continuity of $\uw(\cdot,t_0)$ and $\ow(\cdot,t_0)$ in $(0,L)$,
  by continuity and positivity of $\uw_s(\cdot,t_0)$ and $\ow_s(\cdot,t_0)$ in $(0,L)$, and by local boundedness
  of $\ow_{ss}(\cdot,t_0)$ in $(0,L)\setminus N$ asserted by (\ref{70.5}), we can find $\delta>0$ such that
  \bas
	& & s_j\in [\delta,L-\delta],
	\quad
	\uw(s_j,t_0)\in \Big[-\frac{1}{\delta},\frac{1}{\delta}\Big],
	\quad
	\ow(s_j,t_0)\in \Big[-\frac{1}{\delta},\frac{1}{\delta}\Big], \\
	& &
	\uw_s(s_j,t_0)\in \Big[\delta,\frac{1}{\delta}\Big],
	\quad
	\ow_s(s_j,t_0)\in \Big[\delta,\frac{1}{\delta}\Big]
	\quad \mbox{and} \quad
	\ow_{ss}(s_j,t_0) \in \Big[-\frac{1}{\delta},\frac{1}{\delta}\Big]
  \eas
  for all $j\in\N$.
  As a consequence of this and (\ref{70.3}), there exists $c_2>0$ fulfilling
  \be{70.20}
	|\eta_j| \le c_2
	\qquad \mbox{for all } j\in\N.
  \ee
  Moreover, combining (\ref{70.10}) with (\ref{70.2}), by definition of $c_1$ we obtain that
  \be{70.21}
	|\xi_j|\le c_1
	\qquad \mbox{for all } j\in\N.
  \ee
  Collecting (\ref{70.19}), (\ref{70.20}) and (\ref{70.21}), in (\ref{70.18}) we can further estimate
  \bas
	z_t(s_j,t_0)
	\le c_1 \cdot | \uw(s_j,t_0)-\ow(s_j,t_0)|
	+c_2 \cdot | \uw_s(s_j,t_0)-\ow_s(s_j,t_0)|
	-2c_1 \eps e^{2c_1 t_0}
	\qquad \mbox{for all } j\in\N.
  \eas
  Thanks to the fact that both $\uw$ and $\ow$ belong to $C^1((0,L) \times (0,T_0))$, we may take $j\to\infty$ here
  to see that
  \bas
	z_t(s_0,t_0)
	\le c_1 \cdot | \uw(s_0,t_0)-\ow(s_0,t_0)|
	+c_2 \cdot | \uw_s(s_0,t_0)-\ow_s(s_0,t_0)|
	-2c_1 \eps e^{2c_1 t_0}.
  \eas
  Now observing that $\uw(s_0,t_0)-\ow(s_0,t_0)=\eps e^{2c_1 t_0}$ by (\ref{70.13}), and that
  $\uw_s(s_0,t_0)-\ow_s(s_0,t_0)=0$ by (\ref{70.15}), as a consequence of (\ref{70.14}) we infer that
  \bas
	0\le z_t(s_0,t_0)
	\le c_1 \eps e^{2c_1 t_0}
	-2 c_1 \eps e^{2c_1 t_0}
	< 0.
  \eas
  This absurd conclusion shows that actually (\ref{70.12}) indeed holds, so that on letting $\eps\searrow 0$ and then
  $T_0\nearrow T$ we end up with (\ref{70.9}).
\qed

\end{document}